\documentclass[times,final]{elsarticle}
\usepackage{framed,multirow}

\usepackage[fleqn]{amsmath}
\usepackage{amsthm}
\usepackage{amssymb}
\usepackage{latexsym}

\usepackage{algorithm}
\usepackage{algorithmic}

\usepackage{soul}
\usepackage{lineno}
\usepackage{hhline}

\usepackage{url}
\usepackage{xcolor}
\definecolor{newcolor}{rgb}{.8,.349,.1}

\journal{Elsevier}

\begin{document}

\begin{frontmatter}

\title{Automated tuning for the parameters of linear solvers}

\author{Andrey Petrushov}
\author{Boris Krasnopolsky\corref{cor1}}
\cortext[cor1] {Corresponding author.}
\ead{krasnopolsky@imec.msu.ru}
 
\address{Institute of Mechanics, Lomonosov Moscow State University, Michurinsky ave.~1, Moscow, 119192, Russia}

%\received{\ldots}
%\finalform{\ldots}
%\accepted{\ldots}
%\availableonline{\ldots}

\begin{abstract}
Robust iterative methods for solving large sparse systems of linear algebraic equations often suffer from the problem of optimizing the corresponding tuning parameters. To improve the performance of the problem of interest, specific parameter tuning is required, which in practice can be a time-consuming and tedious task. This paper proposes an optimization algorithm for tuning the numerical method parameters. The algorithm combines the evolution strategy with the pre-trained neural network used to filter the individuals when constructing the new generation. The proposed coupling of two optimization approaches allows to integrate the adaptivity properties of the evolution strategy with a priori knowledge realized by the neural network. The use of the neural network as a preliminary filter allows for significant weakening of the prediction accuracy requirements and reusing the pre-trained network with a wide range of linear systems.

The detailed algorithm efficiency evaluation is performed for a set of model linear systems, including the ones from the SuiteSparse Matrix Collection and the systems from the turbulent flow simulations. The obtained results show that the pre-trained neural network can be effectively reused to optimize parameters for various linear systems, and a significant speedup in the calculations can be achieved at the cost of about 100~trial solves. The hybrid evolution strategy decreases the calculation time by more than 6~times for the black box matrices from the SuiteSparse Matrix Collection and by a factor of 1.4--2 for the sequence of linear systems when modeling turbulent flows. This results in a speedup of up to 1.8~times for the turbulent flow simulations performed in the paper.
\end{abstract}


\end{frontmatter}

%\linenumbers

%%%%%%%%%%%%%%%%%%%%%%%%%%%%%%%%%%%%%%%%

\section{Introduction}
\label{sec:Intro}

The solution of systems of linear algebraic equations (SLAEs) plays an important role in mathematical modeling. Despite the variety of numerical methods developed to date, no single method provides the best possible performance in all cases. The proper selection of the method for the specific problem still remains akin to an art and requires the corresponding qualification of the researcher.

The problem of choosing linear solvers (and solver parameters) has attracted sustained interest by researchers for a long time due to its high practical impact. Several basic research directions dealing with automated optimization problems can be highlighted in the literature. The most popular ones include (\textit{i})~tuning the algorithm implementation details for the specific compute platform, (\textit{ii})~predicting the best numerical method for the specific SLAE, and (\textit{iii})~optimizing the solver configuration for the sequence of SLAEs (e.g., for the transient simulations).

The first research topic deals with fine-tuning the software implementation for the basic computational kernels (linear operations with the vectors, matrix-vector multiplications, etc.), which are the key building blocks of the numerical algorithms~\cite{Whaley2001, Vuduc2005, Dinkelbach2022}. This allows automating the corresponding software adaptation and gaining the benefits of the specific hardware problem.

The two other areas of research mostly focus on algorithmic optimizations and only indirectly touch software implementation issues. The development of assistive algorithms and tools to simplify the choice of efficient solver configurations based on some specific matrix heuristics has been a popular research topic for decades. The tremendous growth of practical interest was associated with developing the machine learning algorithms, which were realized in a series of publications~\cite{Bhowmick2006, Kuefler2008, George2008, Jessup2016, Yeom2016}, to name but a few examples. The corresponding publications investigate various machine learning algorithms, training datasets, lists of features used to characterize the SLAEs, predefined solver configurations, and others. As a result, several software tools (e.g., SALSA~\cite{Demmel2005} and Lighthouse~\cite{Jessup2016}) were developed to recommend the best solver from a predefined set of configurations. These tools, however, being good assistants for novice users, cannot provide highly optimized method configurations with the modified solver parameters set, tuned for the specific SLAE. Moreover, a significant portion of time may be required to calculate matrix heuristics and provide corresponding recommendations.

The third research topic is related to tuning linear solvers (or linear solver parameters) during the calculations when solving a sequence of SLAEs. In~\cite{McInnes2003}, the heuristic algorithm of switching between several predefined linear solver configurations depending on the time integration step was discussed. The same idea is often realized based on the convergence rate history of solving SLAEs at several previous time steps. The machine learning algorithms were applied to predict optimal linear solvers in the transient simulations in~\cite{Eller2012}. An alternative approach based on a genetic algorithm for dynamic solver adaptation was investigated in~\cite{Mishev2008}. The algorithm was used for tuning 3--4~solver parameters with several basic numerical method configurations.

The use of dynamic optimization algorithms is an obvious way to significantly accelerate complex simulations; however, it imposes a number of restrictions when applied to transient problems. The key issue is related to the reproducibility of the simulation results. Some random factors, like minor fluctuations in the calculation time for the specific SLAE, may affect the choice of the method for the next time step and result  in variance in the obtained solution (still satisfying the prescribed solution tolerance). This is especially important when modeling stochastic problems, as any random perturbation eventually brings about a completely different instantaneous solution.

The detailed linear solver parameter optimization for the specific SLAE is a less popular~\cite{GeorgePhD} but still important research topic. These results, for example, can be useful when performing incompressible turbulent flow simulations. The use of high-fidelity models like large eddy simulation or direct numerical simulation (DNS) requires performing long time integration to obtain reliable statistically averaged results. The key computational issue with these models is related to solving the pressure Poisson equation (PPE), which takes up the most of the calculation time. The combination of Krylov subspace and multigrid iterative methods is a popular choice for solving SLAEs corresponding to PPEs. These methods, however, have lots of tuning parameters. For example, the basic set of tuning parameters affecting the multigrid method productivity is counted in several tens, while the per-level multigrid hierarchy parameter specification increases this number to hundreds. The method has no default universal parameters able to provide an acceptable performance level and be applicable for a wide range of SLAEs. The set of parameters optimal for one system can demonstrate poor convergence or divergence for the other SLAEs, even the ones corresponding to the same type of differential equation. Having only rough formal guidelines on how to choose the method parameters, minimization of the SLAE solution time becomes a challenging issue where the researchers can only rely on their experience and intuition. In practice, this issue complicates the large-scale simulations and increases the calculation time, as in most cases, far from optimal linear solver configurations are used to perform the calculations.

The current paper deals with the problem mentioned above and proposes an automated linear solver parameter optimization algorithm. The paper focuses on tuning multigrid-related method parameters and solving SLAEs resulting from the discretization of partial differential equations of the elliptic type with a specific application on modeling incompressible turbulent flows. The method proposed, however, has no strict limitations on the scope of applicability and can be applied to a wider range of applications and areas of mathematical physics.

The rest of the paper is organized as follows. The second section formulates the basic algorithm requirements and limitations caused by the specific problem statement of interest. The proposed optimization algorithm is presented in the third section. The software used for algorithm evaluation and some implementation details, together with the performance evaluation methodology, are summarized in the fourth section. The fifth section presents methodological study results and analyzes the influence of optimization algorithm control parameters. The sixth and seventh sections deal with the proposed algorithm efficiency evaluation. The sixth section investigates the potential of reusing the pre-trained neural networks when optimizing the solver configurations for various linear systems, and the seventh section demonstrates the effect of automated parameter tuning when performing turbulent flow simulations. Finally, the conclusion summarizes the paper.

%%%%%%%%%%%%%%%%%%%%%%%%%%%%%%%%%%%%%%%%

\section{Problem formulation}
\label{sec:problem}

\subsection{Basic requirements}

The current paper deals with developing an automated linear solver parameters optimization algorithm for minimizing the time to solve a sequence of SLAEs
\begin{gather}
\textbf{A} \textbf{x}_i = \textbf{b}_i, \, \, \, \, \, \textbf{A} \in \mathbb{R}^{n \times n}, \, \, \, \textbf{x}_i, \textbf{b}_i \in \mathbb{R}^{n} \label{eq:system}
\end{gather}
occurring in turbulent flow simulations. Accounting for the specific problem statement, the following conditions and limitations are formulated:
\begin{itemize}
  \item as a consequence of the stochastic nature of high-fidelity turbulent flow models and the limitation of the simulation results reproducibility, the optimization procedure can be applied once at the beginning of the corresponding simulation;
  \item the basic numerical methods configuration, the list of optimized parameters, and the parameter optimization range are the input parameters of the optimization algorithm; 
  \item the algorithm must be capable of effectively tuning at least several tens of parameters;
  \item the tuned set of parameters must provide an observable calculation speedup to compensate the time spent on optimization;
  \item it is expected that the linear solver configuration determined for the specific right-hand side (RHS) can be effectively used to solve a sequence of SLAEs occurring when performing the turbulent flow simulation;
  \item some a priori information can be used to increase the efficiency and performance of the optimization algorithm; however, this information must be sufficiently versatile (i.e., valid for at least some group of SLAEs).
\end{itemize}

Without loss of generality, the further narration will be focused on the BiCGStab iterative method~\cite{Vorst1992} with the classical algebraic multigrid preconditioner~\cite{Trottenberg}, a popular combination for solving systems of linear algebraic equations derived from elliptic differential equations. Algebraic multigrid is an example of the method, which has many tuning parameters that have a significant impact on the efficiency when solving specific SLAE.

The problem~\eqref{eq:system} focuses on a series of SLAEs with a constant matrix. This allows for the linear solver initialization stage to be performed only once before the start of calculations. Thus, the real simulation speedup for the problems discussed in the paper is determined by speeding up the linear solver solution stage. All the results presented below account for the corresponding SLAE solution time only, ignoring the time spent on solver initialization. The algorithm developed, however, can be used to optimize the solver while accounting for the initialization time by replacing the evolution function with the sum of the solver's initialization and solution stage times.

%%%%%%%%%%%%%%%%%%%%%%%%%%%%%%

\subsection{Numerical method parameters specification}

In practice, the parameters of the numerical method being optimized can be represented in three ways: as a continuous parameter range, as a range of allowed discrete values, or as a list of allowed discrete values. The proposed optimization algorithm handles all of them in a unified manner. An intermediate parameter mapping is introduced to represent all the parameters in the form of a range of discrete values. The continuous parameter range is discretized with some user-defined stepping to match the corresponding mapping. This approach, however, does not allow for fine-tuning of the parameters between the given points. Estimating the practical significance of this issue as well as possible solutions will be the topic of further research.

%%%%%%%%%%%%%%%%%%%%%%%%%%%%%%%%%%%%%%%%

\section{Optimization algorithm}

Let us denote the following notation. The vector with $N$ parameters is written in the form 
\begin{gather}
\textbf{p} = (p_1, p_2, ... , p_N), \, \textbf{p} \in \mathcal{P}, \\
\mathcal{P} = \{ p_1^i, i=\overline{1,K_1} \} \times \{ p_2^i, i=\overline{1,K_2} \} \times \ldots \times \{ p_N^i, i=\overline{1,K_N} \},
\end{gather}
where $\mathcal{P}$ is the parameter search space, $p_j$ corresponds to a specific numerical method parameter, and each parameter $p_j$ has $K_j$ possible discrete values. The total number of combinations for the parameter vector \textbf{p} equals:
\begin{gather}
K = K_1 \cdot K_2 \cdot ... \cdot K_N.\label{eq:combinations}
\end{gather}
The SLAE solution time with the specific parameter vector \textbf{p} is denoted as
\begin{gather}
T({\textbf{p}}) = T(p_1, p_2, ... , p_N), \label{eq:regression_equation}
\end{gather}
and the role of the optimization algorithm is to find the vector \textbf{p} minimizing the corresponding SLAE solution time.

The form of the expression~\eqref{eq:combinations} clearly shows that simple random search algorithms~\cite{Sivanandam2008} are inapplicable for optimizing problems with several tens of parameters and that more flexible algorithms must be constructed. An example of such an algorithm successfully used for the optimization of the linear solver parameters~\cite{Mishev2008} is the genetic algorithm, a popular variation of the evolutionary algorithms~\cite{Sivanandam2008, Sloss2020}. The method was applied for optimizing problems with 3 and 4~parameters. The current paper, however, deals with a significantly larger number of optimized parameters. Another issue is related to optimization prior to simulation: the test runs to perform the optimization produce overhead for the entire simulation, and the number of these trials must be minimized. 

This paper proposes the hybrid variant of the evolution strategy~(ES)~\cite{Beyer2002}, another popular variation of the evolutionary algorithms. The hybrid ES combines the next-generation vectors produced by the mutation operators with the ones proposed by the neural network~(NN). The basic idea behind this approach is to combine two features: adaptivity to the specific problem of interest provided by the evolution strategy, and reuse of a priori knowledge realized by the neural network.

The algorithm discussed in detail below has an important feature: the proposed scenario of using NN does not require an exact prediction of the specific SLAE solution time, which in practice can be a challenging issue. The key requirement for the NN is only a correct ranking of the parameter configurations. This aspect increases the versatility of the algorithm developed and raises the potential for reusing the pre-trained NN with various SLAEs and compute platforms.

%%%%%%%%%%%%%%%%%%%%%%%%%%

\subsection{Evolution strategy}
\label{sec:ES}
An evolution strategy is a population-based search algorithm suitable for a wide range of optimization problems. A typical ES cycle starts with creating the new generation, which contains several individuals, using the pool of parents. The best individuals are chosen based on the fitness value evaluation and put into the pool of parents to produce the next generation of individuals. The corresponding iterative procedure of constructing the next generation, evaluating fitness values, and selecting the candidates is repeated until some convergence criterion is achieved.

\begin{figure}[t!]
\center{\includegraphics[width=0.65 \textwidth]{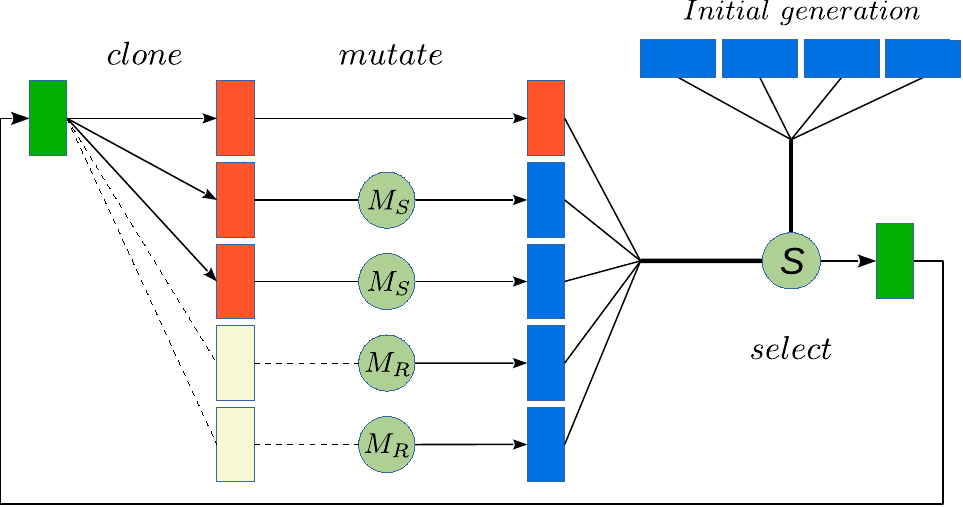}}
\caption{Sketch of the proposed optimization algorithm.}
\label{fig:sketch}
\end{figure}

Following the standard notation~\cite{Beyer2002}, evolution strategy can be denoted as~$(\mu/\rho \overset{\scriptstyle{+}}{,} \lambda)-ES$, where $\mu$ is the size of the parent population, $\rho$ is the mixing number, i.e., the number of parents involved in the formation of new offspring, $\lambda$ is the generation size, and ``+'' or ``,'' indicate the individual selection type. The proposed optimization algorithm is the $(1+\lambda)-ES$ version of the evolution strategy with two mutation operators (Figure~\ref{fig:sketch}). The algorithm uses the following configuration of the ES:
\begin{itemize}
  \item the single parent is used for each individual, $\rho = 1$ (i.e., the new offsprings are clones of their parent);
  \item the soft mutation operator, $M_S$, produces random perturbations in the cloned offspring:
\begin{gather}
\textbf{p} = M_S(\tilde{\textbf{p}}) = ( \tilde{p}_1+s_1, \tilde{p}_2+s_2, ... ,\tilde{p}_N+s_N ),
\end{gather}
where $\tilde{\textbf{p}}$ is the cloned offspring, $s_j \in \{0,\pm \delta_j\}$, and $\delta_j$ is the $j$-th parameter change quantum;
  \item the random mutation operator, $M_R$, is used to introduce some new randomly generated individuals that are generally independent of the cloned offspring;
  \item the selection operator chooses the best individual with the lowest fitness value, $\mu = 1$;
  \item the simple \textit{truncation} technique is applied for choosing the best individual;
  \item the parent is included into the selection process (\textit{plus} selection).
\end{itemize}

The combination of two mutation operators pursues the goal of providing mechanisms for both local adaptation of the current best individual and searching for completely different individuals in other areas of the parameter search space. While playing an important role in overcoming the local minimum problem  (or, at least, finding the one that is closer to the global minimum), the generation of a completely random new individual has significant drawbacks. Some parameter combinations may be incompatible with the linear solver, i.e., the solver would not be able to find the solution of the SLAE with the specified set of numerical method parameters. To avoid falling into ``ineffective'' areas in the parameter search space, the preliminary filtering of the newly generated random input vectors before calculating the fitness value (solving the SLAE) is introduced. The role of the filter is delegated to the pre-trained NN.

%%%%%%%%%%%%%%%%%%%%%%%%%%

\subsection{Neural network model}
\label{sec:NN}

To bring some a priori knowledge about the typical influence of each parameter, a neural network model is built and trained to predict the SLAE solution time for a given vector of parameters. The architecture of the NN is presented in~\ref{app:NN}. The key issues in preparing the input data and estimating the accuracy of the NN predictions are discussed below.

\subsubsection{Training data} 
The SLAE solution time statistic is used to train the NN. In the ideal case, the corresponding data represents a combination of the input parameter vectors and the resulting solution times. In practice, however, some parameter vectors may produce too high SLAE solution times or even lead to solver divergence. The use of NN as a filter with the ES softens the requirements for prediction accuracy: the practical use case scenario requires accurate predictions for the configurations with small execution times only (e.g., values below the median), while for the big ones, rough estimates providing the correct order of magnitude can be acceptable (for details, see Section~\ref{sec:HES}).

The specific problem requirements allow performing initial data filtering in order to improve the predictions accuracy and reduce the training dataset size. The input data preparation includes two modifications: the finite solution times exceeding the third quartile~($Q_3$) and the values for non-converged runs (infinite solution times) are replaced with the $Q_3$ constant. The following notation is used in the further narration: the original dataset is called the ``raw dataset'', the dataset with filtered non-converged runs is called the ``unbalanced dataset'', and the dataset with both filtered non-converged runs and times exceeding the third quartile is called the ``balanced dataset''.

\begin{figure}[t]
\begin{center}
\includegraphics[width=0.7\linewidth]{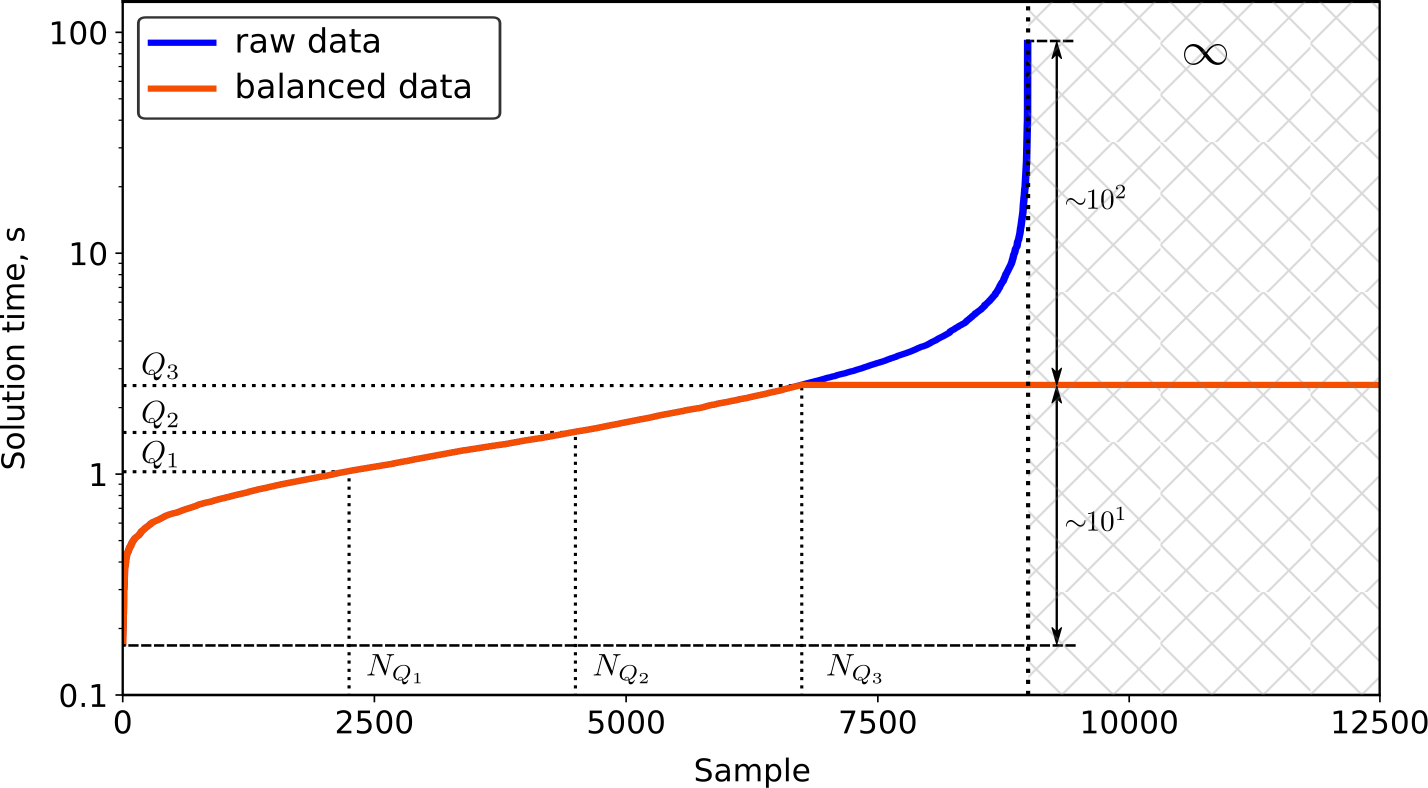}
\end{center}
\vspace{-0.5cm}
\caption{Training dataset initial filtering.}
\label{fig:data_filtering}
\end{figure}

An example of the proposed data filtering is shown in Figure~\ref{fig:data_filtering}. This filtering addresses two issues: (\textit{i})~the formally infinite values for the non-converged runs are replaced with the finite ones, and (\textit{ii})~the dataset is balanced to minimize the prediction errors in the practically important value range. For example, the data modification shown in Figure~\ref{fig:data_filtering} allows for a decrease in the range of predicted values by two orders of magnitude.

\subsubsection{Prediction quality criteria} 
The coefficient of determination, $R^2$, which represents the proportion of the variance explained by the model, is widely used to evaluate the accuracy of the predictive models~\cite{Ahangar2010, Wang2021prediction}. The coefficient is defined as:
\begin{equation}
R^2 = 1-\frac{SS_{R}}{SS_{T}},
\end{equation}
where the residual sum of squares, $SS_{R}$, and the total sum of squares, $SS_{T}$, are:
\begin{gather}
SS_{R} = \sum_{i=1}^{N_v} \left(T^*_i - \hat{T_i}\right)^2,\\
SS_{T} = \sum_{i=1}^{N_v} \left(T^*_i - \bar{T}\right)^2.
\end{gather}
Here $\bar{T}$ is the mean value over the validation dataset
\begin{equation}
\bar{T} = \frac{1}{N_v} \sum_{i=1}^{N_v} T^*_i,
\end{equation}
and $N_v$ is the corresponding dataset size. In cases where the predicted values show a good match with the actual ones, the coefficient of determination tends to be~1. $R^2$ is an integral parameter indicating the proximity of the entire validation dataset and corresponding predictions. While indicating the quality of predictions in general, this metric does not account for the specific requirements of the optimization algorithm being developed, where the key one is the accuracy of predictions in the small value dataset range. 

\begin{figure}[t]
\center{\includegraphics[width=0.6\textwidth]{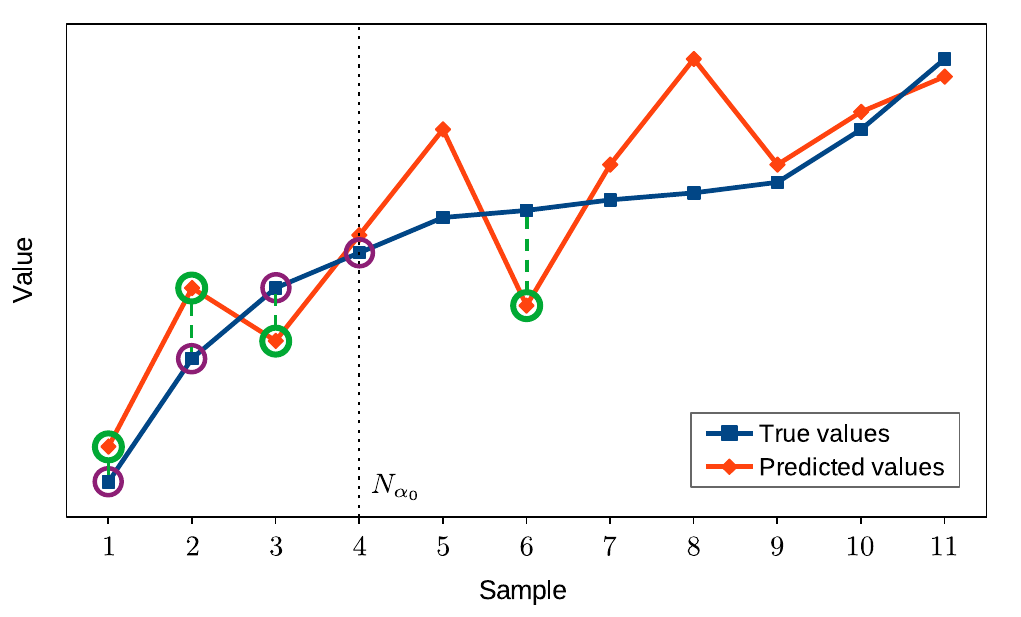}}
\vspace{-0.5cm}
\caption{Prediction quality criteria illustration: $\alpha_0=4/11$, green circles -- least predicted values, purple circles -- least true values; three least predicted values correspond to the least true values.}
\label{fig:prediction_probability}
\end{figure}

An additional customized prediction quality criterion accounting for NN usage specifics is proposed. This metric indicates the probability of finding the least predicted values among the least true ones. The parameter is defined as
\begin{equation}
\mathcal{F}_{\alpha} = \frac{N^{*}_{\alpha}}{N_{\alpha}},
\end{equation}
where $\alpha$ is the fraction of the validation dataset considered, $N_{\alpha} = \alpha N_v$ is the number of samples in the corresponding fraction of the validation dataset, and $N^{*}_{\alpha}$ is the number of least predicted values falling in the list of least true values. The proposed parameter is illustrated in Figure~\ref{fig:prediction_probability}. The plot shows the distribution of some true values sorted in ascending order and the corresponding predictions. According to this figure, for $\alpha_0=4/11$ the size of the least values list is~$N_{\alpha_0}=4$, the number of samples both falling in the lists of least predicted and true values is equal to~$N^{*}_{\alpha_0}=3$, and $\mathcal{F}_{\alpha_0} = 0.75$.

%%%%%%%%%%%%%%%%%%%%%%%%%%%%%%%%%%%%%%%%%%%%%%%%

\subsection{Hybrid optimization algorithm}
\label{sec:HES}
The proposed hybrid ES (HES) combines both the ES and NN optimization techniques discussed above. The basic $(1+\lambda)-ES$ optimization algorithm is presented in Algorithm~\ref{alg:HES}, and the random mutation operator incorporating the NN is outlined in Algorithm~\ref{alg:M_R}. 

% The random mutation operator, if configured to use the neural network, generates a large preliminary set of parameter 
The random mutation operator with neural network pre-filter generates a large preliminary set of parameter configurations containing $L$~samples. The predictions for these samples are obtained using the NN, and $L_{\alpha} = \alpha L$ samples with the least predicted values are selected. $\lambda_R$ vectors are randomly chosen from the selected list and passed to the ES.

In the limiting case, $L_{\alpha}$ equals to $\lambda_R$, and $L = \lambda_R / \alpha$ provides the lower bound estimate for the trial dataset size. This, however, is not recommended as it can affect the probability of finding the least true values in the least predicted values list (see Section~\ref{sec:F_alpha}). To produce reliable results for the $\mathcal{F}_{\alpha}$ criteria, the size of the preliminary set of combinations should be chosen in accordance with the size of the validation dataset.

The random selection of $\lambda_R$ vectors from the larger least value parameters list presumes the same probability of finding the true least values, and $\lambda_R$ vectors would contain on average $\mathcal{F}_{\alpha} \lambda_R$ parameter combinations from the true least values list in each ES generation. At least one true least value is expected to be passed to the ES in each generation. The reasonable upper bound of the trial dataset can be obtained based on the typical validation dataset size required for an accurate enough prediction of the parameter $\mathcal{F}_{\alpha}$. The corresponding question is discussed in detail in Section~\ref{sec:F_alpha}.

The algorithm formulated above operates with the parameter $\alpha$. In general, selecting $\alpha$ as little as possible can be preferable, as it will shrink the predicted least values list and preserve the parameters with the smallest predicted values. However, the decrease of $\alpha$ also affects the $\mathcal{F}_{\alpha}$ parameter, and the probability of predicting true least values decreases. The optimal value $\alpha^*$ can be found by maximizing the ratio
\begin{equation}
P(\alpha) = \frac{\mathcal{F}_{\alpha}}{\alpha}\label{eq:alpha_estimate}
\end{equation}
with $\mathcal{F}_{\alpha} \geq 1/\lambda_R$ as a constraint. In practice, a moderate value, such as $\alpha^* = 0.05$, can be used as a baseline value, and it can be tuned in each case depending on the quality of predictions produced by the specific NN.

\begin{algorithm}[t]
\caption{Hybrid evolution strategy $(1+\lambda)-ES$.} 
\begin{algorithmic}[1]
\label{alg:HES}
\STATE $l = 0$;
\STATE $P^0 = \{ \textbf{p}_1^0, \textbf{p}_2^0, \ldots,  \textbf{p}_{\lambda}^0 \}$ -- initial generation with $\lambda$ individuals;
\STATE Evaluate fitness values for $P^0$ and select the best individual $\textbf{p}^{0,*}$;
    \WHILE {(\textit{convergence criteria is not met})}
    \STATE Create $\lambda$ offsprings $\tilde{\textbf{p}}^{l+1}$ by cloning the selected individual $\textbf{p}^{l,*}$
    \STATE Apply soft mutation operator, $M_S$, to $\lambda_S$ individuals: \\ $\textbf{p}^{l+1}_k = M_S(\tilde{\textbf{p}}^{l+1}_k), \, \, k=\overline{1, \lambda_S}$;
    \STATE Apply random mutation operator, $M_R$, to $\lambda_R$ individuals (Algorithm~\ref{alg:M_R}): \\ 
           $\textbf{p}^{l+1}_k = M_R(\tilde{\textbf{p}}^{l+1}_k), \, \, k=\overline{\lambda_S + 1, \lambda_S+\lambda_R}$, \\ where $\lambda_S+\lambda_R = \lambda$;
    \STATE Evaluate fitness values for the new generation $P^{l+1}$ and select the best individual $\textbf{p}^{l+1,*}$ among $\textbf{p}^{l+1}_k$ and $\textbf{p}^{l,*}$;
    \STATE $l = l+1$;
\ENDWHILE 
\end{algorithmic}
\end{algorithm}

\begin{algorithm}[t]
\caption{Random mutation operator, $M_R$.} 
\begin{algorithmic}[1]
\label{alg:M_R}
\IF{(\textit{use NN filtering})}
  \STATE Generate $L$ trial vectors $\tilde{\textbf{p}}_k, \, k=\overline{1,L}$, $L \gg \lambda_R$;
  \STATE Predict solution ``times'' $T(\tilde{\textbf{p}}_k), \, k=\overline{1,L}$ with the NN;
  \STATE Choose $L_{\alpha}$ vectors with the least predicted values;\label{alg:choose}
  \STATE Randomly choose $\lambda_R$ vectors from $L_{\alpha}$ candidates selected at the step~\ref{alg:choose}, assign vectors as new individuals in ES;
\ELSE
\STATE Generate $\lambda_R$ vector parameters $\tilde{\textbf{p}}_k, \, k=\overline{1,\lambda_R}$;
\STATE Assign generated parameters as new individuals in ES.
\ENDIF
\end{algorithmic}
\end{algorithm}

To complete, the ES termination criteria must be provided. A couple of options can be used in practice, e.g., non-decreasing or slowly decreasing fitness value during the last $m$~iterations. The current paper uses the condition of non-decreasing fitness value in the last 5~iterations together with the limitation on the overall number of generations (50~generations).

%%%%%%%%%%%%%%%%%%%%%%%%%%%%%%%%%%%%%%%%%%%%%%%%%%%%

\section{Algorithm evaluation details}

The proposed hybrid evolution strategy has been tested in detail with a set of model SLAEs and integrated with in-house code for modeling incompressible turbulent flows. The software implementation details and the test matrices used for testing are presented below.

\subsection{Implementation details}

The HES is implemented as a part of the XAMG library~\cite{xamg_repo} for solving SLAEs with multiple right-hand sides. The library provides implementation of a set of iterative methods, including the Krylov subspace and algebraic multigrid methods, and supports mixed precision calculations~\cite{KrasnopolskyIR2021}. The library is C++ code based on the C++11 standard specification. It provides hierarchical multilevel parallelization with the MPI+POSIX shared memory hybrid programming model. Some additional XAMG library implementation details and performance evaluation results can be found, e.g., in~\cite{Krasnopolsky2021xamg, Krasnopolsky2021rsd}.

The TensorFlow framework~\cite{tensorflow2015} is used to operate with NN. The NN training process is organized as a standalone task, performed separately with the Python code developed. The inference stage is incorporated into the XAMG library using the TensorFlow C API, and the corresponding functionality is directly used from the XAMG library.

\subsection{Compute platforms}

The calculations presented in the paper have been performed on cluster compute nodes with 2 x 14-core Intel Xeon E5-2680~v4 processors and on the workstation with an 8-core Intel Core i7-11700~processor. Unless explicitly stated, the test results are obtained with a single cluster compute node. The hybrid 3-level solver configuration corresponding to the hardware architecture for both compute platforms is used in the XAMG library. All available CPU cores per node are utilized in the calculations.

\subsection{Test SLAEs}
\label{sect:test_matrices}
A set of linear systems predominantly corresponding to the elliptic differential equations is used to evaluate the efficiency of the proposed optimization algorithm. The test set includes several model SLAEs widely used to evaluate the efficiency of numerical methods (e.g.,~\cite{Yang2010, Gahvari2013, Bienz2016}), a subset of matrices from the SuiteSparse Matrix Collection~\cite{SuiteSparse}, and a group of SLAEs coming from the turbulent flow simulations~\cite{KrasnopolskyCPC2018}.

% The first group of linear systems includes two SLAEs corresponding to 3D problems in a cubic computational domain. Regular grid and 7-point discretization stencil are used in both cases to solve the diffusion equation
% \begin{gather}
% - \nabla \cdot \left( \mathcal{K} \nabla u \right) = f, \label{eqn:diffusion}
% \end{gather} 
% with Dirichlet boundary conditions, and the cases considered differ in the form of the diffusion coefficient, $\mathcal{K}$. 

The first group of linear systems includes two SLAEs corresponding to 3D problems in a cubic computational domain. The diffusion equation
\begin{gather}
- \nabla \cdot \left( \mathcal{K} \nabla u \right) = f, \label{eqn:diffusion}
\end{gather} 
with Dirichlet boundary conditions is discretized on the regular grid using the 7-point stencil, and the cases considered differ in the form of the diffusion coefficient, $\mathcal{K}$. The first case denoted as \texttt{cube} corresponds to the constant diffusion coefficient, $\mathcal{K} = 1$. The second case, \texttt{jumps}, considers the stepwise distribution of the diffusion coefficient:
\begin{gather}
\mathcal{K}(\textbf{x}) = 
\begin{cases}
1000, & \textbf{x} \in [0.1, 0.9]^3, \\
0.1, & \textbf{x} \in \mbox{cubes of size} \, 0.1^3 \, \mbox{at the corners of the domain}, \\
1, & \mbox{elsewhere}.
\end{cases}
\end{gather}

The second group includes a subset of the matrices from the SuiteSparse Matrix Collection. Square matrices with nonzero diagonal elements and the number of unknowns in the range of $10^5$ to $2\cdot10^6$ are selected and supplemented with the constant right-hand side vector $\textbf{b} = \textbf{1}$. The ones that can be solved with the \textit{baseline} BiCGStab method and algebraic multigrid preconditioner configuration (for details, see~\ref{app:list}) are included in the SSMC test set for HES evaluation.

Finally, the third group of test SLAEs relates to the systems occurring when modeling incompressible turbulent flows, the target problem of the current research. The SLAEs correspond to the solution of PPE arising in direct numerical simulation of turbulent flow in a channel with a matrix of wall-mounted cubes~\cite{KrasnopolskyCPC2018}. The two cases, \texttt{DNS1} and \texttt{DNS2}, correspond to two structured hexahedral grids of size 2.3 and 9.7~million unknowns, respectively.

\subsection{Numerical method configurations}

The three sets containing 7, 13, and 21 optimized parameters are considered in most of the tests for evaluating the efficiency of the HES. An additional 15-parameter configuration is used for methodological studies in Section~\ref{sec:param_significance}. The parameters are related to the classical algebraic multigrid and Chebyshev polynomial methods used as a preconditioner and smoother with the BiCGStab solver, respectively. The lists of parameters included in these optimized sets are shown in \ref{app:list}.

%%%%%%%%%%%%%%%%%%%%%%%%%%%%%%%%%%%%%%%%

\section{Methodological studies}

The detailed methodological studies are performed to evaluate the efficiency of the proposed HES and the typical range of control parameters affecting its productivity and performance. These include both the parameters of the evolution strategy and the neural network-related parameters. These questions are discussed in detail in the following sections.

The test results presented below deal with the NNs trained with the data for three test systems: two \texttt{cube} systems constructed for $40^3$ and $100^3$~grid points, and the \texttt{apache2} system from the SSMC. The three linear solver configurations with 7, 13, and 21 parameters are considered for each system. The randomly generated parameter vectors are used to produce SLAE solution time statistics and to train the neural network model. In total, 50~thousand samples are collected for each system and parameter set. It should be noted that the obtained datasets sufficiently differ in the number of non-converged samples: for example, for the \texttt{cube40} system with a 7-parameter set, all the combinations produce finite solution times, while for the \texttt{apache2} system with a 13-parameter set, the fraction of non-converged combinations exceeds~66\%.

%%%%%%%%%%%%%%%%%%%%%%%%%%%%%%%%%%%%%%%%

\subsection{Neural network training}

The methodological studies investigate in detail the aspects of neural network training. This section analyzes the influence of the control parameters when training the NN and estimates the accuracy of NN predictions.

\subsubsection{$\mathcal{F}_{\alpha}$ prediction quality criteria}
\label{sec:F_alpha}
The methodological studies start with estimating the prediction quality criteria for the trained neural networks and analyzing the required validation dataset size to obtain reliable $\mathcal{F}_{\alpha}$ estimates. The tests use the NN trained with 45~thousand samples, and the other 5~thousand samples are used for validation. The most complex configuration, with 21~parameters, is considered. To estimate the spread of results obtained, in accordance with the cross-validation procedure, the corresponding learning process and $\mathcal{F}_{\alpha}$ calculations are repeated five times, and the mean over the five trials as well as the standard deviation are presented.

\begin{figure}[t]
\center{\includegraphics[width=\textwidth]{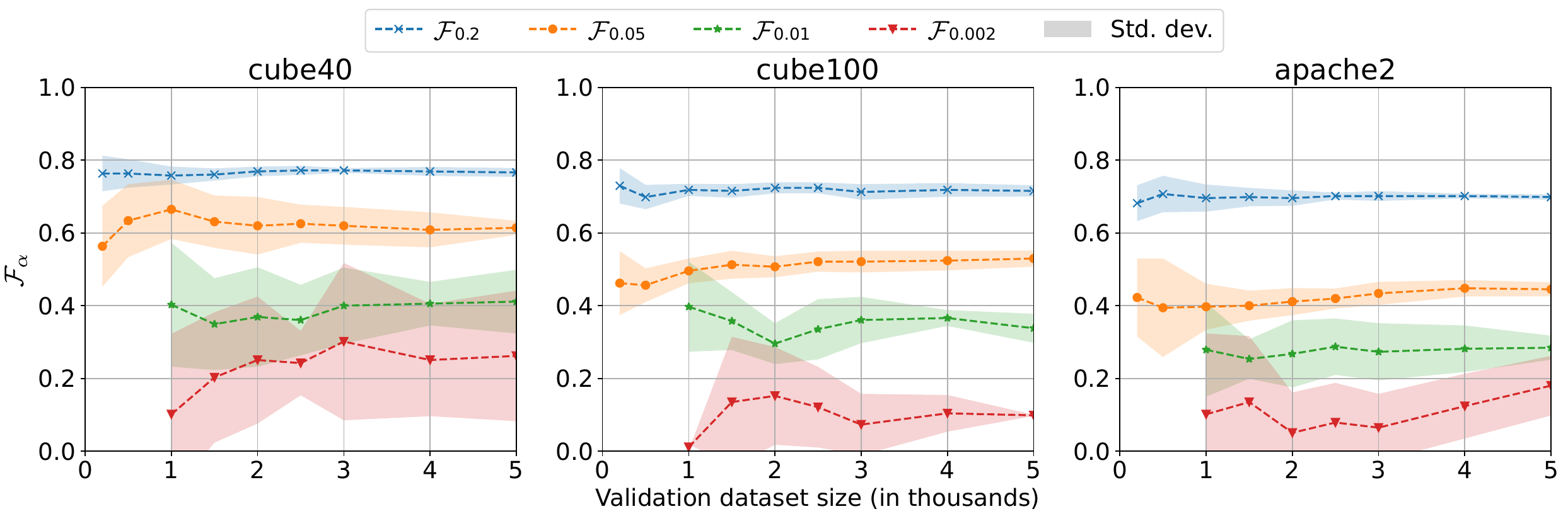}}
\vspace{-0.5cm}
\caption{Dependence of the $\mathcal{F}_{\alpha}$ prediction quality metric on the validation dataset size, $N_v$.}
\label{fig:FL_cone}
\end{figure}

The obtained results for $\mathcal{F}_{\alpha}$ with $\alpha = 0.002$, $0.01$, $0.05$, and $0.2$ as a function of the validation dataset size are shown in Figure~\ref{fig:FL_cone}. The results presented demonstrate that the spread over trials decreases with increasing the validation dataset size, and the dataset with $N_v = 5000$~samples (in most cases, even lower) allows for reliable $\mathcal{F}_{\alpha}$ estimates. As expected, lowering the fraction of the selected least predicted values leads to a decrease in the absolute values of the prediction quality criteria. However, in accordance with~\eqref{eq:alpha_estimate}, the obtained values for $\mathcal{F}_{0.002}$ can still be used in the HES with a large enough number of randomly mutated individuals.

\subsubsection{Data balancing}

\begin{figure}[t]
\center{\includegraphics[width=0.8\textwidth]{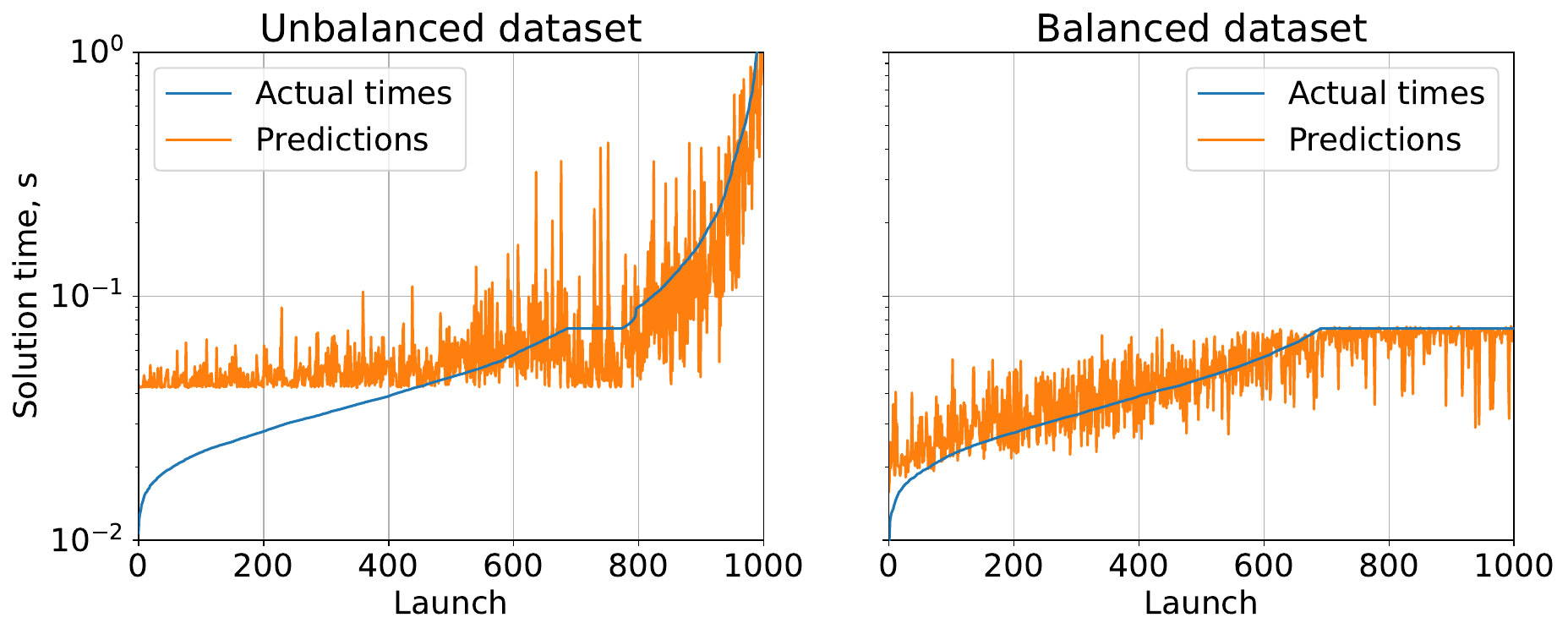}}
\vspace{-0.3cm}
\caption{Comparing the quality of predictions for unbalanced and balanced datasets. Test case: \texttt{cube40}, 21~parameters set.}
\label{fig:balancing_results}
\end{figure}

Data balancing plays an essential role when training the neural network. The specific balancing procedure proposed in Section~\ref{sec:NN} helps improve the accuracy of predictions. The preliminary methodological studies have shown that the presence of the extremely large SLAE solution times in the dataset significantly coarsens the predictions for the small ones. An example of the corresponding training is the one obtained for the unbalanced data containing 20~thousand samples for the configuration with 21~parameters and the \texttt{cube40} test system. Figure~\ref{fig:balancing_results} shows the real solution times from the validation dataset sorted in ascending order and NN predictions. The plot clearly demonstrates the weak quality of the predictions for the unbalanced data and the tendency to predict some constant value. The data filtering procedure applied to the dataset before NN training replaces the longest solution times. This allows to reduce the spread of the solution times (for example, in Figure \ref{fig:balancing_results} -- from 0.01--3~s to 0.01--0.073~s). Training with the balanced dataset significantly improves the accuracy of predicting the lowest values for the validation dataset. The predicted values correctly reproduce the general trend for the practically important combinations when used with HES.

The effect of data balancing is also illustrated with the $\mathcal{F}_{0.05}$ and $R^2$ metrics calculated for three test matrices and three parameter sets with training datasets containing 45~thousand samples (Table~\ref{tab:balancing_table_45K}). Both datasets with unbalanced and balanced data allow for reliable results for 7- and 13-parameter solver configurations. The extended configuration with 21~parameters, on the other hand, clearly shows the data balancing effect: the use of balanced input data allows for a significant increase in both the $\mathcal{F}_{0.05}$ and $R^2$ results. An increase in the training dataset size improves the quality of predictions in all cases, and the most notable one is the 21-parameter solver configuration.

It should be noted that the $\mathcal{F}_{\alpha}$ parameter shows more sensitive behavior with regard to data balancing compared with $R^2$. Using unbalanced input data for the 21-parameter case and the \texttt{cube40} test system, for example, provides a rather good value, $R^2=0.73$, whereas the proposed metric yields only $\mathcal{F}_{0.05} = 0.28$. Balancing the data raises $R^2$ by 13\% while increasing $\mathcal{F}_{0.05}$ by 130\%. This aspect indicates the need to control the $\mathcal{F}_{\alpha}$ metric in addition to $R^2$ when training the NN for use with HES.

\begin{table*}[t]
\caption{Comparing the NN prediction quality criteria for unbalanced and balanced data, 45~thousand samples.}
\centering
\begin{tabular}{| c | c | c | c | c | c | c |}
\hline
Solver & Test & \multicolumn{2}{c|}{Unbalanced}  & \multicolumn{2}{c|}{Balanced} & Improvement, \\
\cline{3-6}
 configuration & system & $\mathcal{F}_{0.05}$  & $R^2$ & $\mathcal{F}_{0.05}$  & $R^2$ & $\mathcal{F}_{0.05}$, in \% \\
\hline
\multirow{3}*{7 parameters} & \texttt{cube40} & 0.84 & 0.99 & 0.87 & 0.99 & 3.9\\
& \texttt{cube100} & 0.85 & 0.99 & 0.86 & 0.98 & 0.8\\
& \texttt{apache2} & 0.86 & 0.97 & 0.88 & 0.98 & 2.4\\
\hline
\multirow{3}*{13 parameters} & \texttt{cube40} & 0.72 & 0.92 & 0.68  & 0.88 & -5.1\\
& \texttt{cube100} & 0.72 & 0.85 & 0.76 & 0.86 & 6.0\\
& \texttt{apache2} & 0.65 & 0.65 & 0.67 & 0.75 & 2.1\\
\hline
\multirow{3}*{21 parameters} & \texttt{cube40} & 0.28 & 0.73 & 0.65 & 0.83 & \textbf{129.3}\\
& \texttt{cube100} & 0.20 & 0.52 & 0.57  & 0.70 & \textbf{189.4}\\
& \texttt{apache2} & 0.19 & 0.40 & 0.49 & 0.58 & \textbf{165.1}\\
\hline
\end{tabular}
\label{tab:balancing_table_45K}
\end{table*}

%%%%%%%%%%%%%%%%%%%%%%%%%%%%%%%%%%%%%%%%%%

\subsubsection{Training dataset size estimates}
The dependence of the characteristic dataset size for training the neural network on the number of optimized parameters is investigated. A series of tests with progressively increasing dataset sizes is performed to estimate the corresponding requirements. The $\mathcal{F}_{\alpha}$ parameter is evaluated in these tests, and the validation dataset with 5~thousand samples is used to calculate $\mathcal{F}_{\alpha}$.

\begin{figure}[t]
\centering
\includegraphics[width=\linewidth]{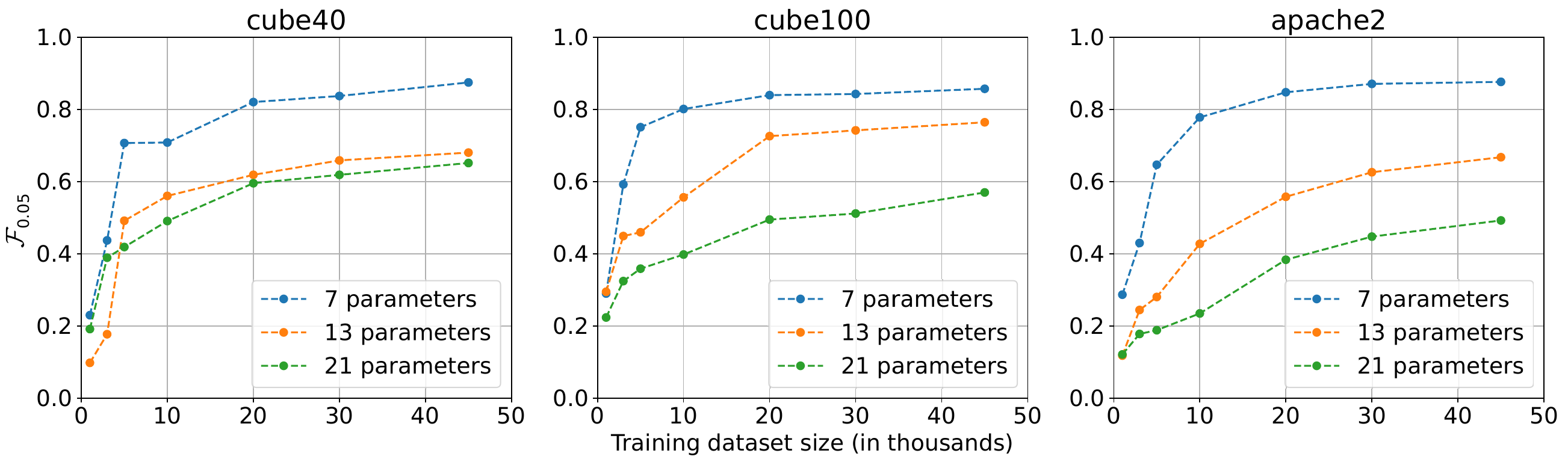}
\vspace{-0.8cm}
\caption{Dependence of the $\mathcal{F}_{0.05}$ prediction quality metric on the training dataset size.}\label{fig:NN_training_dataset}
\end{figure}

Figure~\ref{fig:NN_training_dataset} shows the typical testing results. The plots presented for $\alpha=0.05$ demonstrate acceptable quality of predictions starting from 20~thousand samples, and further increasing the size of the training dataset only slightly increases the corresponding parameter values. A similar picture is observed when varying the $\alpha$ parameter: only a minor difference in the absolute values takes place, while the corresponding trends remain the same. These results indicate that the dataset with 45~thousand samples, used in the tests presented below, is sufficient for training the NN and obtaining reliable predictions.

%%%%%%%%%%%%%%%%%%%%%%%%%%%%%%%%%%%%%%%%%%%%%%%%%%%%%%%%

\subsection{Hybrid evolution strategy}

HES has several control parameters affecting the overall performance. These include the generation size and the number of soft and random mutations used to produce new individuals. In addition to them, the fraction of the least values range must be specified when using NN as a random mutation filter. The following sections investigate the effect of varying these parameters, their impact on the optimization results, and provide some basic recommendations for choosing these values.

The randomization process is an essential part of the evolution strategy. This, however, may result in some variation in the optimization results from run to run. The calculation results presented below were obtained by averaging over 20~trials, and the mean and standard deviation are typically shown.

The results presented in the following subsections are performed for three test systems, \texttt{cube40}, \texttt{cube100}, and \texttt{apache2}, and for three solver configurations with 7, 13, and 21 parameters. The NNs used when optimizing the linear solver configuration are trained with datasets containing information for the same SLAE and for the same compute system.

\subsubsection{Neural network as a pre-filter}
\label{sec:NN_filter}
The first test series demonstrates the effect of using the NN as a pre-filter in the random mutation operator. The tests evaluate the HES with 5~soft and 10~random mutations, and the random mutation operators both with and without NN are considered. For NN-assisted cases, the parameter $\alpha=0.05$ is used. In addition to optimal SLAE solution times, the HES calculation times spent on finding the optimum are also collected and analyzed.

\begin{figure}[t!]
\centering
\includegraphics[width=\textwidth]{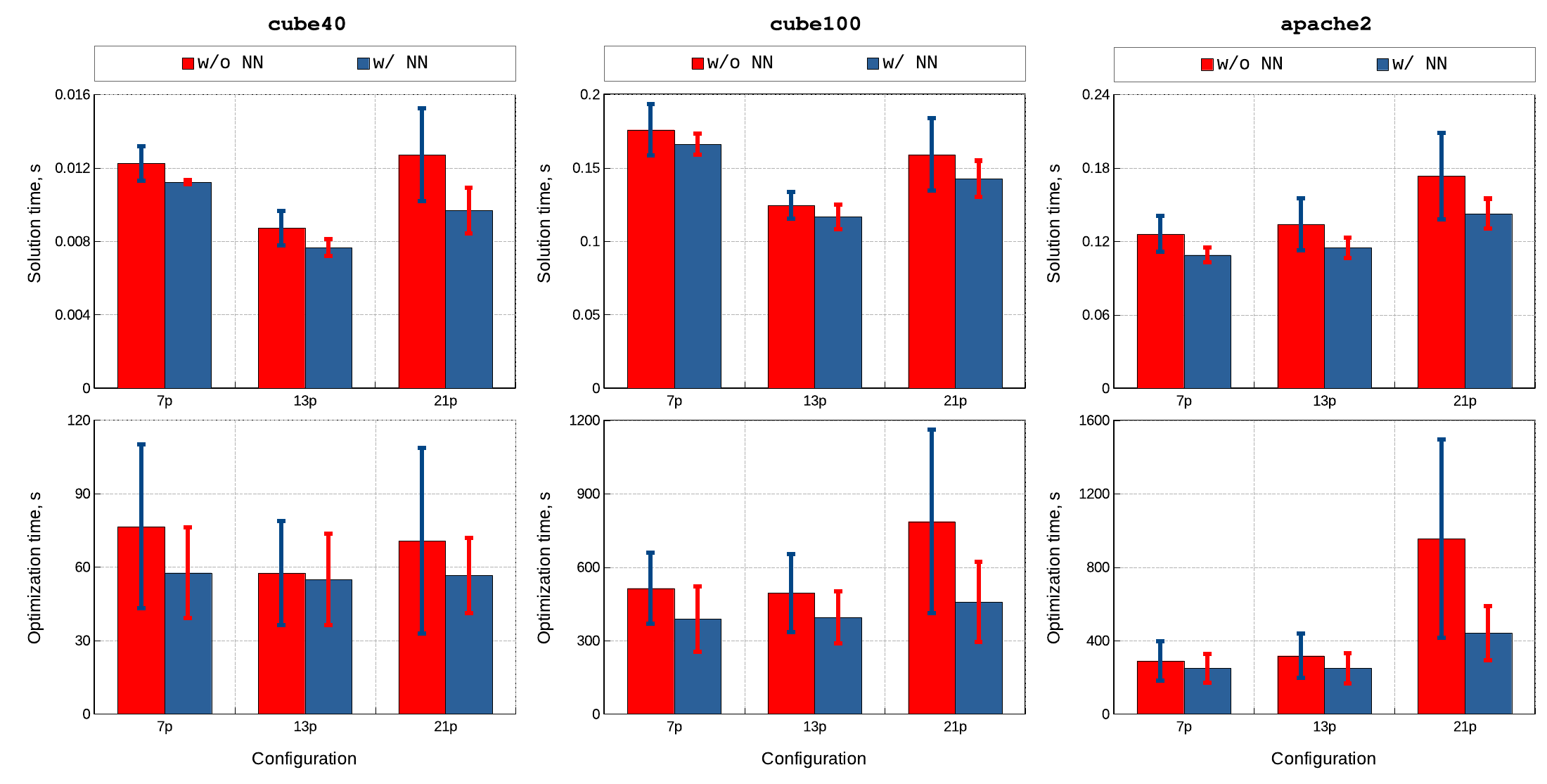}
\vspace{-0.75cm}
\caption{Effect of using NN as a pre-filter in the random mutation operator: top -- optimized solution times, bottom -- HES calculation times.}
\label{fig:NN_filter}
\end{figure}

The obtained results (Figure~\ref{fig:NN_filter}) confirm the initial expectations regarding the use of NN as a pre-filter:
\begin{itemize}
  \item the NN-assisted random mutation operator improves the optimized SLAE solution time; for the SLAEs considered, the solution time reduction reaches 30\%;
  \item the NN filtering allows for a decrease in the HES calculation time, and a speedup by a factor of~1.7 can be achieved;
  \item the spread for both the solution and calculation times is significantly decreased when using the NN filter.
\end{itemize}

The use of NN does not change the mean number of generations of HES until reaching convergence. The convergence for the criteria specified in Section~\ref{sec:HES} is typically achieved in 10~generations, which is equivalent to 100-200~trial SLAE solutions.

%%%%%%%%%%%%%%%%%%%%%%%%%%%%%%%%%%%%%%%

\subsubsection{Soft and random mutations}
The generation size and the number of soft and random mutations are the important parameters affecting the quality of optimization and optimization time. A detailed study is performed to estimate the optimal ratio between the soft and random mutations. Four configurations are considered:
\begin{itemize}
  \item \texttt{S5/R5}: 5 soft and 5 random mutations;
  \item \texttt{S5/R10}: 5 soft and 10 random mutations;
  \item \texttt{S10/R5}: 10 soft and 5 random mutations;
  \item \texttt{S10/R10}: 10 soft and 10 random mutations.
\end{itemize}
The NN filtering is used with the random mutation operator, and the fraction of the least values range is set equal to $\alpha=0.05$.

\begin{figure}[t]
\centering
\includegraphics[width=\textwidth]{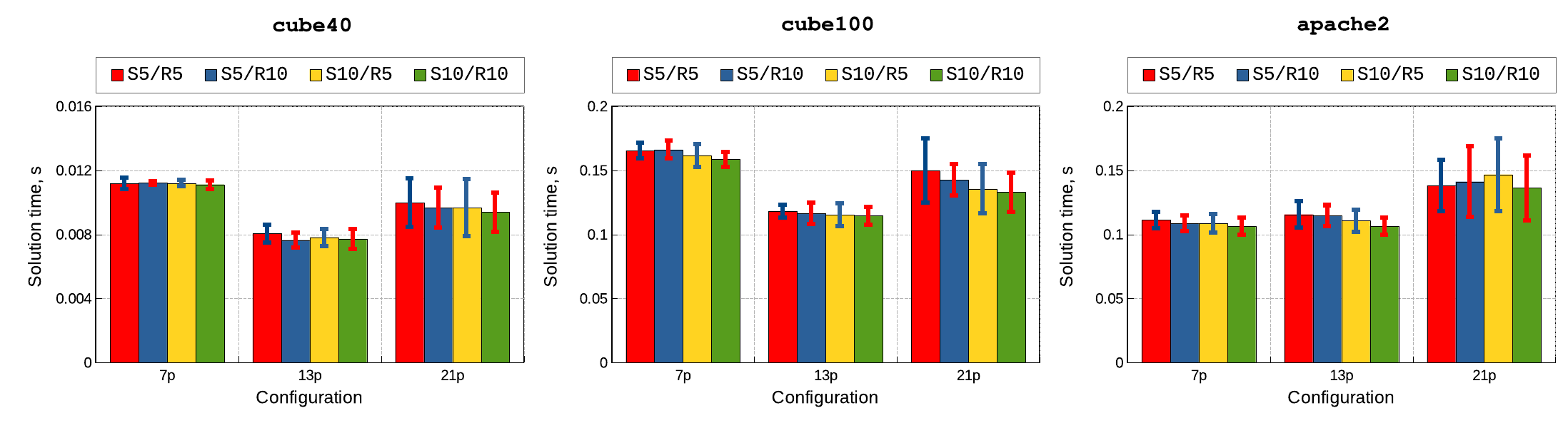}
\vspace{-0.75cm}
\caption{Effect of varying the number of soft and random mutations.}
\label{fig:mutations}
\end{figure}

The analysis of the calculation results shows that an increase in the number of individuals (generation size) reduces the resulting optimized SLAE solution time (Figure~\ref{fig:mutations}). However, any substantial difference in SLAE solution times is observed for the 21-parameter case only when the variance across the test configurations reaches 10\%. In all other cases, the use of different HES configurations leads to negligible variation in the optimized results, but the algorithm calculation time increases in proportion to the generation size. These results allow us to conclude that the use of the \texttt{S5/R5} HES configuration can be an acceptable choice when using the optimization algorithm with the NN-assisted mutation operator for the algorithm efficiency evaluation purposes discussed in this paper.

%%%%%%%%%%%%%%%%%%%%%%%%%%%%%%%%%%%%%%%

\subsubsection{Choosing the fraction of the least values range}
\label{sec:var_alpha}
One more parameter used with NNs is the fraction of the least predicted values when producing the individuals for the NN-assisted random mutation operator. The simple relation for determining the optimal value is proposed in~\eqref{eq:alpha_estimate}. Following the results presented in Section~\ref{sec:F_alpha}, the value of $\alpha^* = 0.002$ must be close to the optimal one for all the NNs trained in this work. The numerical experiments performed focus on validating these predictions.

\begin{figure}[t]
\centering
\includegraphics[width=\textwidth]{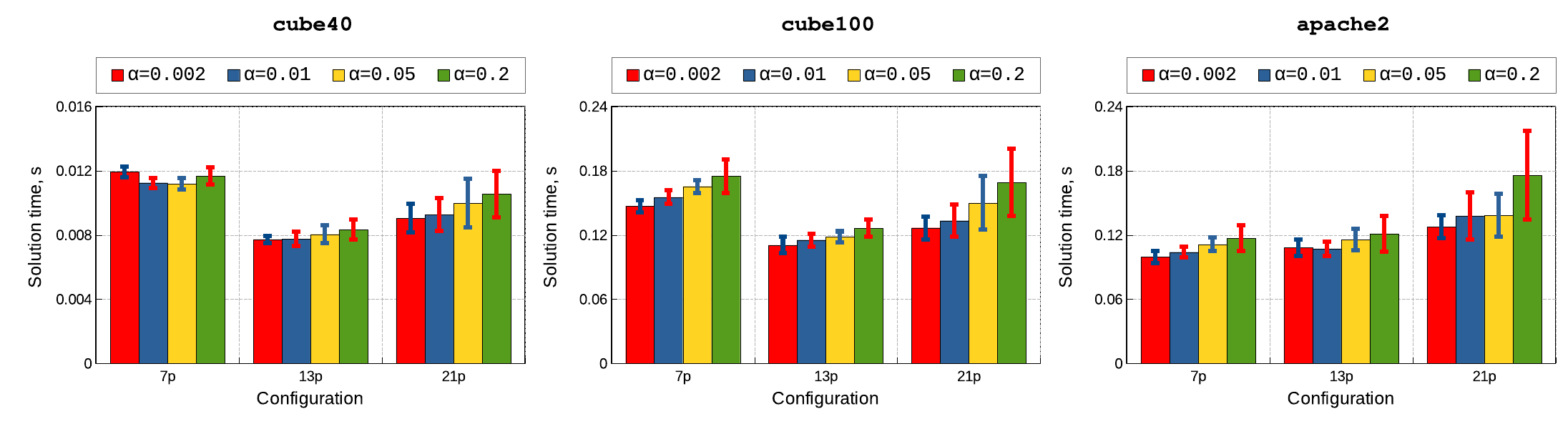}
\vspace{-0.75cm}
\caption{Effect of varying the fraction of least values range, predicted by the NN.}
\label{fig:alpha}
\end{figure}

This test series performs the experiments for the HES configuration with five soft and five random mutations and with four $\alpha$ values: 0.002, 0.01, 0.05, and 0.2. The corresponding results are summarized in Figure~\ref{fig:alpha}. The clear tendency to improve the quality of optimization with decreasing the $\alpha$ parameter can be seen (except only for the \texttt{cube40} system and 7-parameter solver configuration), which coincides with predictions~\eqref{eq:alpha_estimate}. Reducing $\alpha$ from 0.2 to 0.002 results in a 25\% reduction in SLAE solution time while maintaining the HES calculation time.

%%%%%%%%%%%%%%%%%%%%%%%%%%%%%%%%%%%%%%%

\subsection{Convergence of the optimization algorithm}

Analyzing the solution times as a function of the number of optimized parameters, one can expect a reduction in the resulting times with an increase in the parameter list. This takes place when switching from 7 to 13~parameters, but the resulting times for the 21-parameter runs can exceed the times for the 13-parameter ones (Figure~\ref{fig:NN_filter}). An increase in the number of parameters in practice has a dual effect: an extended set of parameters accounts for more features of the numerical method and can provide lower solution times. At the same time, an increased number of parameters hardens the optimization process due to a significant growth in the possible combinations (see Table~\ref{tab:4sets}). Thus, an extension of the list of optimized parameters raises the problem of balancing the potential benefits of increasing the parameter set with some degradation of the quality of optimization. The following subsections briefly discuss this issue.

\subsubsection{Parameters significance}
\label{sec:param_significance}

Decreasing the efficiency of optimization for the 21-parameter set can indicate that the 8~\textit{extra} parameters (which extend the 13-parameter set to~21) already have sub-optimal default values: their variation does not lead to improvement in the solution times but complicates the optimization of the other parameters. This supposition is verified with an additional test series performed for the new 15-parameter set constructed as a combination of the 7-parameter set and 8~\textit{extra} parameters (for details, see~\ref{app:list}). The standard procedure for collecting the data and training the NN is performed for the \texttt{cube100} problem, which allowed to achieve the quality of training $R^2=0.95$ and $\mathcal{F}_{0.05} = 0.77$. The ES configuration and testing conditions correspond to those used in Section~\ref{sec:NN_filter}.

\begin{figure}[t!]
\centering
\includegraphics[width=0.5 \textwidth]{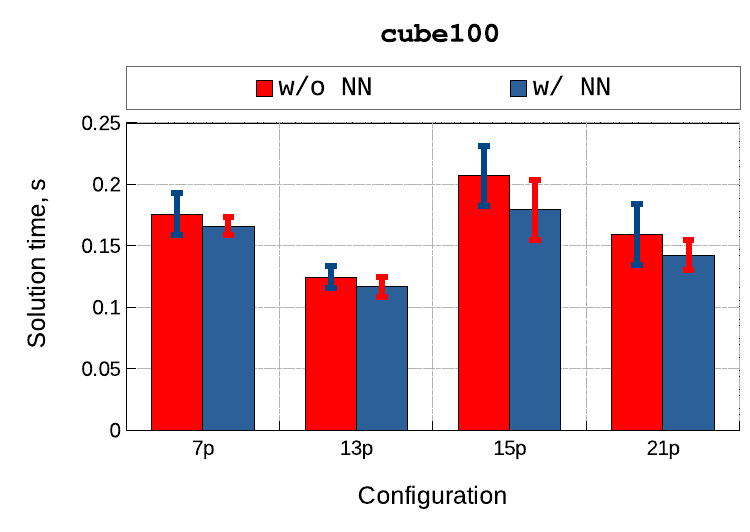}
\vspace{-0.5cm}
\caption{Comparison the optimization results for different sets of optimized parameters.}
\label{fig:set15}
\end{figure}

Testing results shown in Figure~\ref{fig:set15} demonstrate a similar effect when adding 8~\textit{extra} parameters to the sets with 7 and 13 parameters. For both ES and HES, an increase in the solution times as well as in the standard deviations is noticed. Extending the 7-parameter configuration with different parameters shows the opposite trend: adding six parameters allows for an observable reduction in the SLAE solution time, while adding another eight parameters increases the obtained optimized values. These observations comply with the suggestions discussed above and highlight the importance of the problem of parameter classification and identifying the reduced list of parameters, which strongly affect the solution times.

%%%%%%%%%%%%%%%%%%%%%%%%%%%%%%%%%%%%%%%

\subsubsection{Stopping criteria}

\begin{figure}[t!]
\centering
\includegraphics[width=\linewidth]{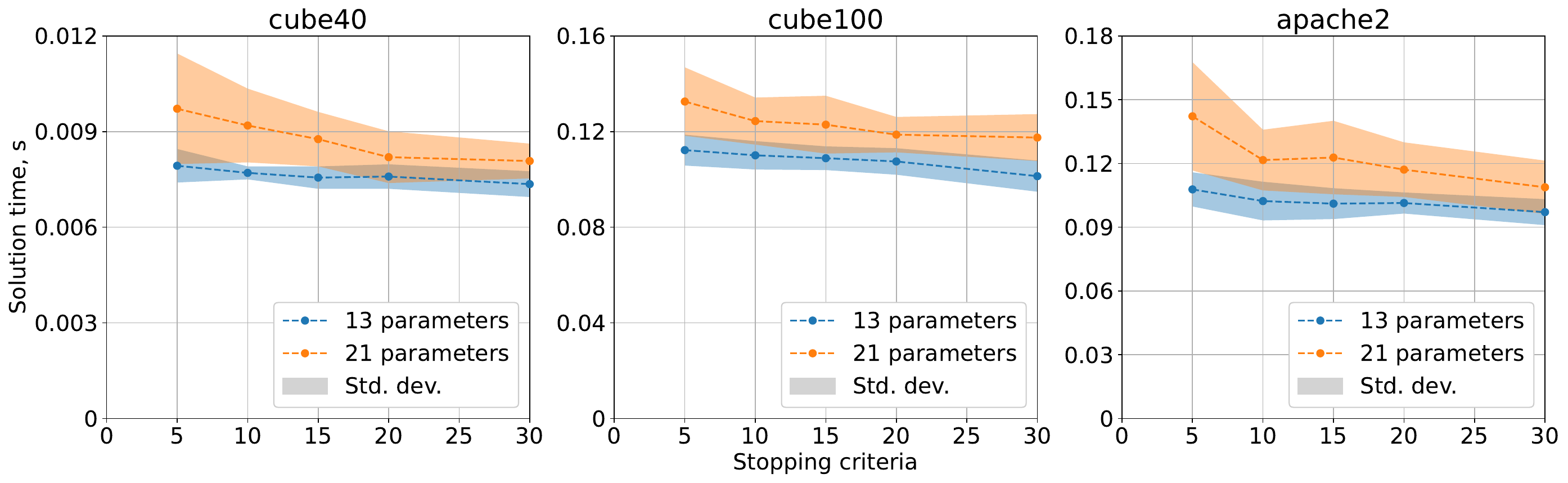}
\vspace{-0.7cm}
\caption{Dependence of the optimized solution times on the number of generations used for the stopping criteria.}
\label{fig:stopping_criteria_triple}
\end{figure}

The observed difference in solution times when adding \textit{extra} parameters is related to the worsening quality of the optimization for the extended parameter set. The two tests are performed to evaluate the effect of varying the convergence criteria. The same configuration of the evolution strategy is used as in Section~\ref{sec:NN_filter}. The first test series performs the calculations with the progressively strengthening termination criteria: starting from 5~generations of non-decreasing optimal solution times, this limit is increased to~30. The optimized solution time, together with the standard deviation for 13 and 21 parameters, are shown in Figure~\ref{fig:stopping_criteria_triple}. Increasing the convergence limit leads to improvements in optimized solution time and a decrease in standard deviation. The most notable improvement is observed for the 21-parameter configuration when switching from 5 to 10~generations. This may indicate that an increase in the number of generations can be reasonable when increasing the list of optimized parameters. Also, the tendency to decrease the deviation between the configurations with 13~and 21~parameters can be seen.

\begin{figure}[t!]
\centering
\includegraphics[width=0.75\textwidth]{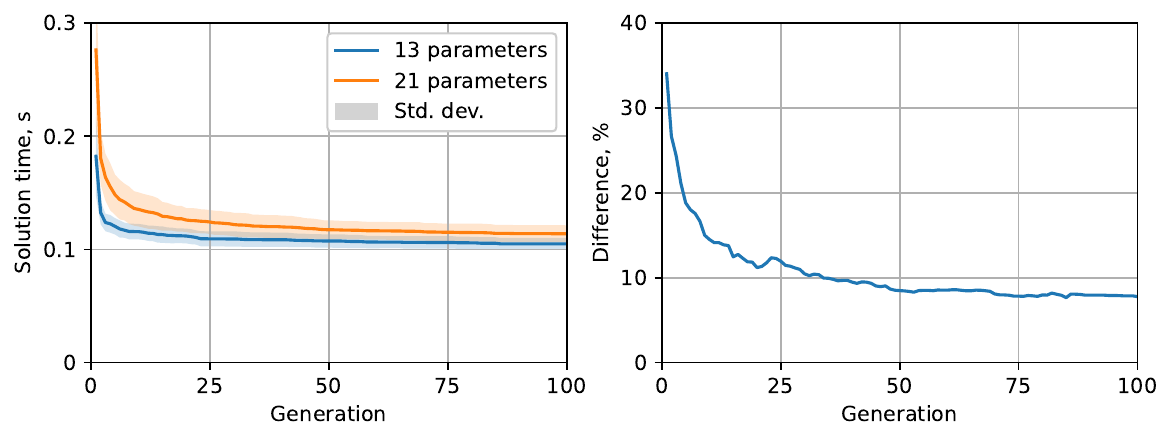}
\vspace{-0.3cm}
\caption{Evolution of the solution times and their difference for the 13- and 21-parameter configurations as a function of the generation number of the HES, \texttt{cube100} case.
\label{fig:conv}}
\end{figure}

The testing scenario above controls the number of generations with a non-changing objective function. An additional test series analyzes the evolution of the mean values as a function of the number of generations for the \texttt{cube100} test case. The corresponding results in Figure~\ref{fig:conv} indicate that the difference between the configurations considered has a tendency to decrease, dropping to about 10\% in 20-30~generations. The further increase in the number of generations until~100 allows for reducing this difference below~8\%. 

%%%%%%%%%%%%%%%%%%%%%%%%%%%%%%%%%%%%%%%

\subsubsection{Additional convergence improvements}

An alternative approach with multistage optimization, realizing an adapted number of individuals per generation or the number of soft and random mutations, parameter change stepping, etc., can be considered. The use of an initial guess for the optimization algorithm to produce the first generation of individuals can also improve the convergence. These optimizations, however, are out of scope for this paper and are the subject of future work.

%%%%%%%%%%%%%%%%%%%%%%%%%%%%%%%%%%%%%%%%%%%%%%%%%%%%%%%%

\section{Reusing the pre-trained neural networks}

The tests shown above have been performed using the NNs trained with the same linear systems statistics. While allowing for the investigation of the influence of the optimization algorithm control parameters, the corresponding configuration is an idealization of the practical use cases. The following section investigates the perspectives of reusing the pre-trained NNs to optimize linear solver configurations for various linear systems.

%%%%%%%%%%%%%%%%%%%%%%%%%%%%%%%%%%%%%%%

\subsection{Test systems cross-optimization}
The first test series performs the cross-optimization for the basic test systems (\texttt{cube40}, \texttt{cube100}, and \texttt{apache2}) and the NNs trained for the same matrices to estimate the potential of reusing the pre-trained NNs. \texttt{Cube40} and \texttt{cube100} systems have the same origin, but differ in the computational grid and matrix size by about 6~times; the \texttt{cube100} and \texttt{apache2} matrices are of about the same size, but produced by various differential equations. The \texttt{S5/R5} HES configuration is used in these tests, and two values, $\alpha=0.002$ and $\alpha=0.05$, are compared.

\begin{figure}[t]
\centering
\includegraphics[width=\textwidth]{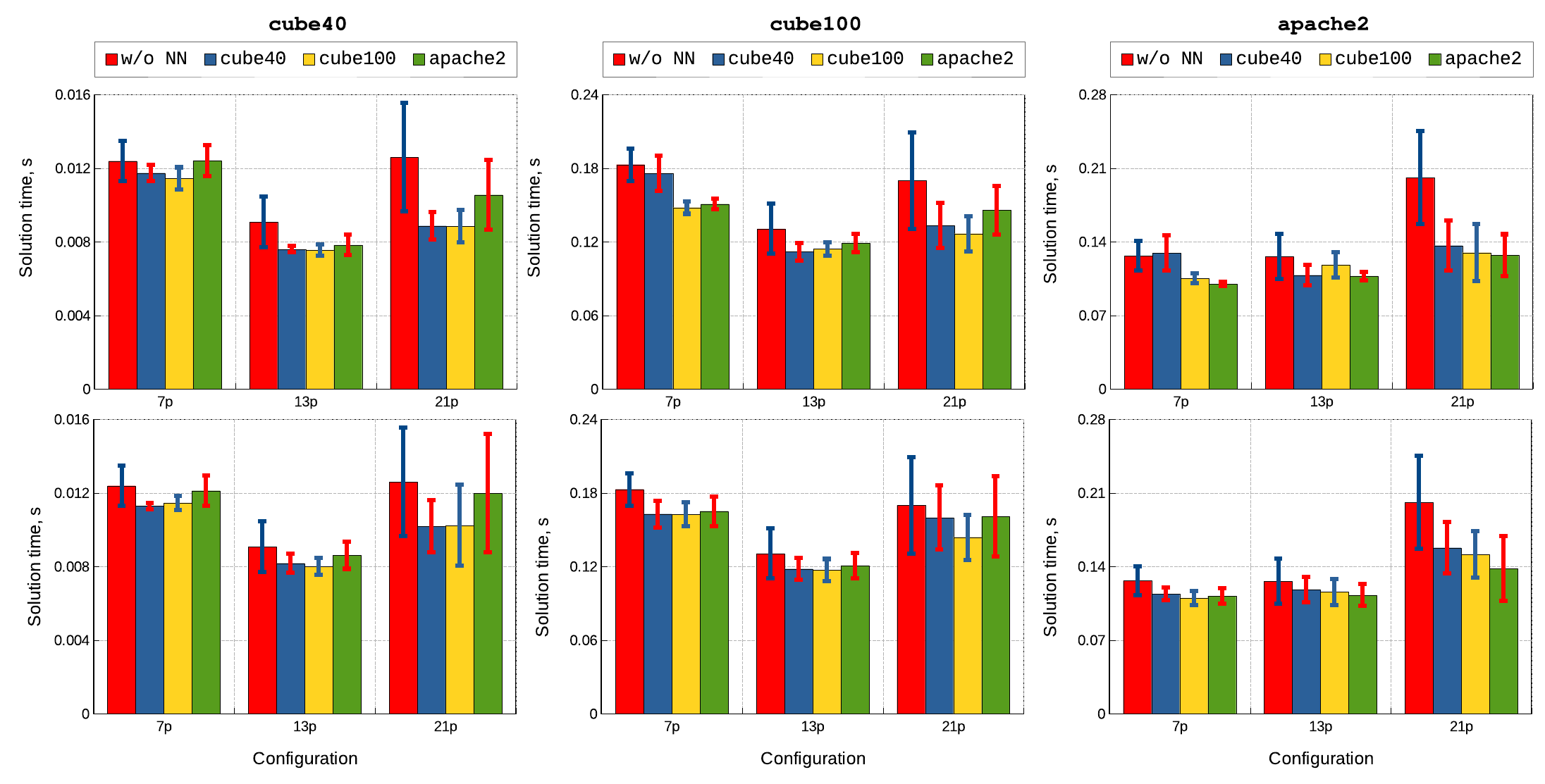}
\vspace{-0.75cm}
\caption{Cross-optimization using various pre-trained NNs: top -- $\alpha = 0.002$, bottom -- $\alpha = 0.05$.}
\label{fig:cross_opt}
\end{figure}

The obtained calculation results confirm an important HES feature -- the possibility of reusing the pre-trained NNs with various linear systems (Figure~\ref{fig:cross_opt}). The optimization of the \texttt{cube40} SLAE solution time can be done with the same effect with both the \texttt{cube40} and \texttt{cube100} pre-trained NNs. The \texttt{apache2} NN demonstrates lower solution time improvement; however, it still outperforms the ES without the NN pre-filter. Additionally, the NN pre-filter with any of the NNs considered significantly decreases the variation across the trial runs.

The \texttt{cube100} SLAE can be successfully optimized with any NN among the considered ones. The native \texttt{cube100} NN demonstrates slightly better results, but the variance of the corresponding times does not exceed several percent. The same situation is observed for the \texttt{apache2} system: the evident advantage of using the native \texttt{apache2} pre-trained NN is seen for the 21-parameter case with $\alpha=0.05$, while all other test configurations show negligible differences in the optimized solution times.

A comparison of the results for various fractions of the least values range reveals that using $\alpha=0.002$ is superior to $\alpha=0.05$ and allows for a 5\% reduction in solution time as well as variation across trials. This observation is consistent with the results discussed in Section~\ref{sec:var_alpha}.

%%%%%%%%%%%%%%%%%%%%%%%%%%%%%%%%%%%%%%%

\subsection{Linear systems for elliptic differential equations}
\label{sec:elliptic}
The second test series performs the optimization of SLAE solution times for the \texttt{jumps}, \texttt{DNS1}, and \texttt{DNS2} systems. These systems correspond to elliptic differential equations and are similar to those that are of interest to us. The tests also use the \texttt{S5/R5} HES configuration with $\alpha=0.002$ and two pre-trained NNs for the \texttt{cube100} and \texttt{apache2} systems. These calculations are done with 4~compute nodes. In addition to optimized SLAE solution times, the results for two predefined \textit{baseline} and \textit{default} linear solver configurations formed on the basis of some calculation experience and \textit{hypre} library recommendations are evaluated.

\begin{figure}[t]
\centering
\includegraphics[width=\textwidth]{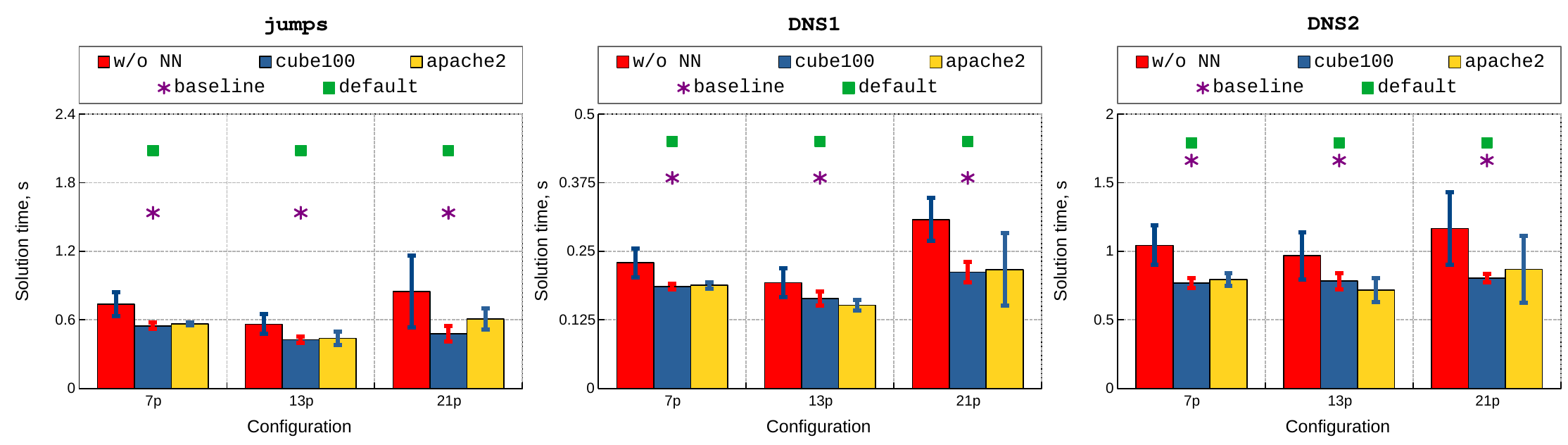}
\vspace{-0.75cm}
\caption{Linear solver configurations optimization results for various test systems.}
\label{fig:elliptic_optimization}
\end{figure}

The obtained calculation results are summarized in Figure~\ref{fig:elliptic_optimization}. These results show that the provided \textit{baseline} parameter configuration demonstrates systematically better results compared with the \textit{default} one. The proposed evolution strategy can significantly improve \textit{baseline} times and tailor the linear solver configuration to the specific linear system. The use of \texttt{cube100} or \texttt{apache2} pre-trained NNs improves the optimization results even further. The choice of the optimal NN varies depending on the linear system and optimized parameter set, but the variance of the mean values is within 10-15~percent.  Despite this, the use of any of the NNs considered allows for a reduced optimized solution time compared to the ES. The speedup by a factor of 1.3--2.7 is observed for the ES, and it increases to 1.8--3.6 compared to the \textit{baseline} when using the pre-trained NN. In addition, the use of the NN pre-filter for the random mutation operator reduces the algorithm calculation times, and this effect becomes more noticeable with increasing the number of optimized parameters. Up to a 5-fold HES calculation time speedup when optimizing the 21-parameter configuration is observed with the \texttt{cube100} NN compared to a basic ES run.

% The best solution times after the optimization are observed for the 13-parameter configuration. The use of a more general 21-parameter case, which further extends the parameter search space, typically does not allow for any productivity improvement. This fact indicates that some further optimization algorithm tuning or more conscious choice of the optimizing parameters may produce additional productivity gains.

%%%%%%%%%%%%%%%%%%%%%%%%%%%%%%%%%%%%%%%
	
\subsection{Reusing NNs across various platforms}
\label{sec:NN_reuse}\label{Ed:summary1}
The tests above were performed with the NNs constructed with the data collected on the same hardware platform. The current test series evaluates the possibility of reusing the pre-trained NNs across multiple compute platforms. The tests below operate with two NNs constructed for the 13-parameter set, using the training data for the \texttt{cube100} case on the cluster node (2~processors, 24~cores) and on the workstation (8~cores).

\begin{figure}[t]
\centering
\includegraphics[width=7.0cm]{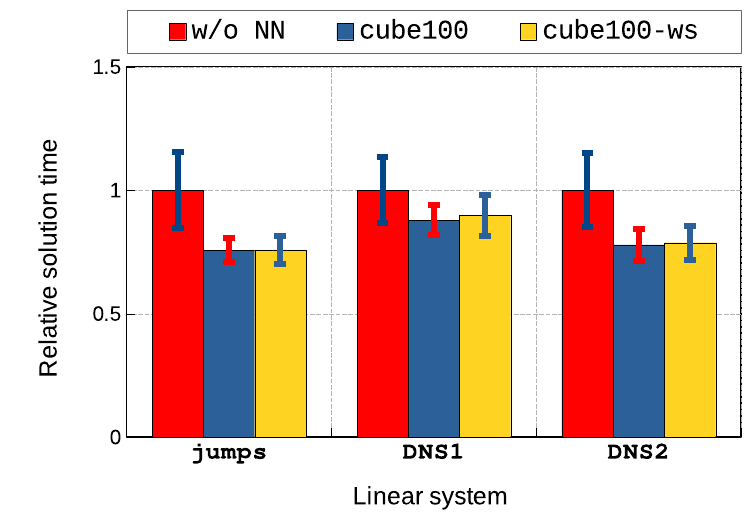}
\vspace{-0.5cm}
\caption{Relative optimized solution times for various test problems when using NNs trained with the data for different compute platforms: \texttt{cube100} -- NN trained with the cluster node data, \texttt{cube100-ws} -- NN trained with the workstation data.}
\label{fig:NN_ws}
\end{figure}

The new test series is performed for \texttt{jumps}, \texttt{DNS1}, and \texttt{DNS2} problems on four cluster nodes to compare the optimized solution times for the ES and the HES with two NNs. The corresponding results, presented in the form of relative solution times normalized with the ES optimized solution times, are shown in Figure~\ref{fig:NN_ws}. The obtained results demonstrate negligible variation in the resulting mean solution times within the statistical averaging accuracy. These results give reasons to believe that the pre-trained NNs can be effectively reused across various compute platforms (at least, similar ones in some sense). This can also indicate that the key optimization effect when optimizing the linear solver configuration is related to algorithmic speedup and is meaningful across multiple compute platforms but not the software implementation details, which are compute platform-specific.

Minor variation in the solution times can be seen if comparing the results in Figures~\ref{fig:elliptic_optimization} and~\ref{fig:NN_ws}. This is due to performing different test runs. The variation of the mean values and the standard deviations across them is within 1-2\%, which complies with the expected accuracy of the results shown in the paper.

%%%%%%%%%%%%%%%%%%%%%%%%%%%%%%%%%%%%%%%

\subsection{SuiteSparse Matrix Collection}

The final test series is performed for the SSMC systems. Unlike the cases considered above, these SLAEs come from different areas of applied mathematics and do not necessarily correspond to elliptic differential equations. These tests exceed the main purpose of the developed algorithm and evaluate its ability to work with arbitrary linear systems. An attempt to optimize the solution time with the HES is made for the selected SLAEs, and the \texttt{cube100} and \texttt{apache2} pre-trained NNs are evaluated. Following the methodological study results, the \texttt{S5/R5} HES configuration is chosen, and the fraction of the least predicted values for the hybrid random mutation operator is set to~0.002.

\begin{figure}[t]
\centering
\includegraphics[width=11cm]{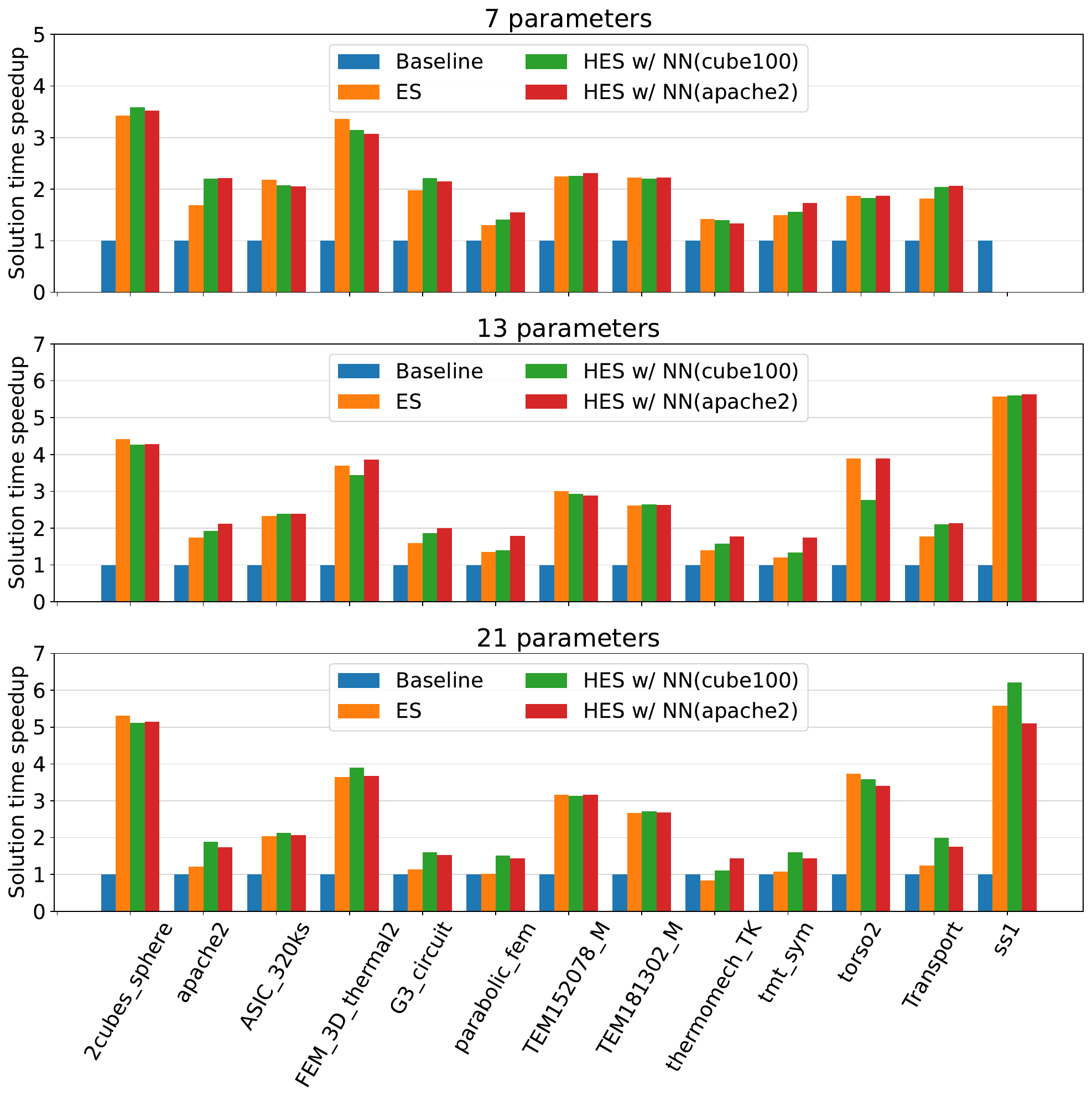}
\caption{Optimizing the linear solver configurations for the SSMC test matrices; missing data for the 7-parameter configuration and the \textit{ss1} system correspond to the non-converging ES and HES cases.}
\label{fig:SSMC_results}
\end{figure}

The numerical experiments conducted have shown that the use of the evolution strategy allows for a significant improvement in the SLAE solution time (Figure~\ref{fig:SSMC_results}). In all cases except the \texttt{ss1} matrix with 7-parameter set, HES produces the solver parameter configurations capable of solving the test SLAEs. In the worst case (\texttt{thermomech\_TK}, 21~parameters, ES), the obtained solution time exceeds the \textit{baseline} by 16\%. In all other cases, the \textit{baseline} times are significantly improved, and the peak acceleration by a factor of 6.2 is observed for the \texttt{ss1} system with a 21-parameter configuration. The average SLAE solution time reduction of the ES for the configurations with 7, 13, and 21~parameters reaches 2.18, 2.42, and~2.26, respectively.

The use of NNs provides limited advantages in terms of optimized solution time and calculation time for a set of black box matrices. The minor SLAE solution time improvement by a factor of 1.05--1.2 can be seen. The same behavior takes place for the optimization time: both the slowdown by a factor of~1.1 and the speedup by a factor of~1.4 are observed for different test configurations. The NN filtering typically gives no advantage when the ES shows significant solution time improvement. However, it allows for optimizing the SLAE solver configurations for the cases when the ES fails to produce anything better than the \textit{baseline} results. On the one hand, this observation indicates the importance of using NN filtering with the HES. On the other hand, a more careful selection of the pre-trained NNs may be required to further improve the ES productivity when working with black box linear systems. The obtained results comply with initial expectations -- the optimal multigrid method parameters vary significantly for the problems of various types, which makes it almost impossible to predict optimal ones using the NN trained with the single linear system data. The possible ways of combining multiple NNs or training the NN with multiple SLAE results are planned to be investigated in the future.

%%%%%%%%%%%%%%%%%%%%%%%%%%%%%%%%%%%%%%%%%%%%%%%%%%%%%%%%

\section{Direct numerical simulation of turbulent flow}

The second part of the performance evaluation tests considers the modeling of incompressible turbulent flow. The direct numerical simulation of the turbulent flow is performed, and the potential calculation speedup due to automated optimization of the linear solver parameters is demonstrated. The problems simulated are the modeling of incompressible turbulent flow in the plain channel, turbulent flow in a channel with a matrix of wall-mounted cubes~\cite{KrasnopolskyCPC2018, Meinders1999}, and turbulent flow in a straight duct~\cite{Gavrilakis1992}. The 13-parameter solver configuration is used for the simulations following the methodological tests performed as it is the most suitable in terms of resulting SLAE solution time, optimization time, and the spread of results from run to run. The evolution strategy uses the~\texttt{S5/R5} configuration for the first two cases with the \texttt{cube100} pre-trained NN, and $\alpha$ is set equal to~0.002. The last case uses an increased number of individuals and an~\texttt{S10/R10} algorithm configuration.

The in-house computational code is used to perform DNS. The code is based on a second-order spatial finite difference scheme operating with curvilinear orthogonal coordinates and a third-order Runge-Kutta scheme for advancing in time~\cite{Nikitin2006jcp, Nikitin2006ijnmf}. The time integration scheme used requires three solves of the Poisson equation per time step, which takes up the overwhelming majority of the calculations. The XAMG library is used to solve these linear systems. 

The two runs are performed for each test case with different linear solver configurations: the \textit{baseline} configuration and the one produced by the HES optimization algorithm. In addition, ES is used to perform the calculation of the turbulent flow in a square duct. The corresponding calculation times and the number of linear solver iterations are analyzed to demonstrate the possibility of using the optimized solver configurations to solve a series of SLAEs. The evolution of the iteration number in the figures is correlated with the friction at the walls to indicate the changes when performing the transition to a statistically stationary turbulent flow regime. Here, the non-dimensional skin friction coefficient is defined as
\begin{gather}
C_f = \frac{\tau_w}{0.5 \rho U_b^2},
\end{gather}
where $\tau_w$ is the wall shear stress, $\mbox{Re}_b = U_b h / \nu$ is the bulk Reynolds number, $U_b$ is the bulk velocity, $h$ is the characteristic length scale, and $\nu$ is the kinematic viscosity. 

All the times mentioned below are represented in non-dimensional global time units (scaled by $h / U_b$). The specific time integration intervals for the tests considered are chosen for demonstration purposes only. They correlate by order of magnitude but are somewhat lower than those typically used in the DNS to obtain reliable statistical averaging~\cite{Vinuesa2016}.

%%%%%%%%%%%%%%%%%%%%%%%%%%%%%%%%%%%%%%%%

\subsection{Plain channel flow}
\label{sec:channel}
The direct numerical simulation of the turbulent flow in a plain channel with $\mbox{Re}_b = 2800$ ($\mbox{Re}_{\tau} = 180$) is considered. The flow is modeled in a $2\pi h \times 2h \times \pi h$ box with a grid of $160 \times 140 \times 160$ cells, where $h$ is the channel half-height. The grid is uniform in streamwise and spanwise directions and stretched in the wall-normal direction. The problem statement is similar to the one used in~\cite{KrasnopolskyCPC2018}.

The simulation starts with perturbed laminar Poiseuille flow; the evolution of these perturbations leads to a transition to turbulence. About 200~time units must be modeled to obtain statistically stationary turbulent flow, and further integration is performed to collect statistics. In these tests, the overall duration of the simulation time interval is set equal to $T=500$ time units with a time integration step of $\tau = 0.02$. The tolerances for solving the SLAEs (both absolute and relative residual, whatever occurs first) are set equal to $10^{-6}$ for the pressure and $10^{-8}$ for the pseudo-pressure~\cite{Nikitin2006jcp}. These calculations are performed on five compute nodes.

\begin{figure}[t]
\centering
\includegraphics[width=11cm]{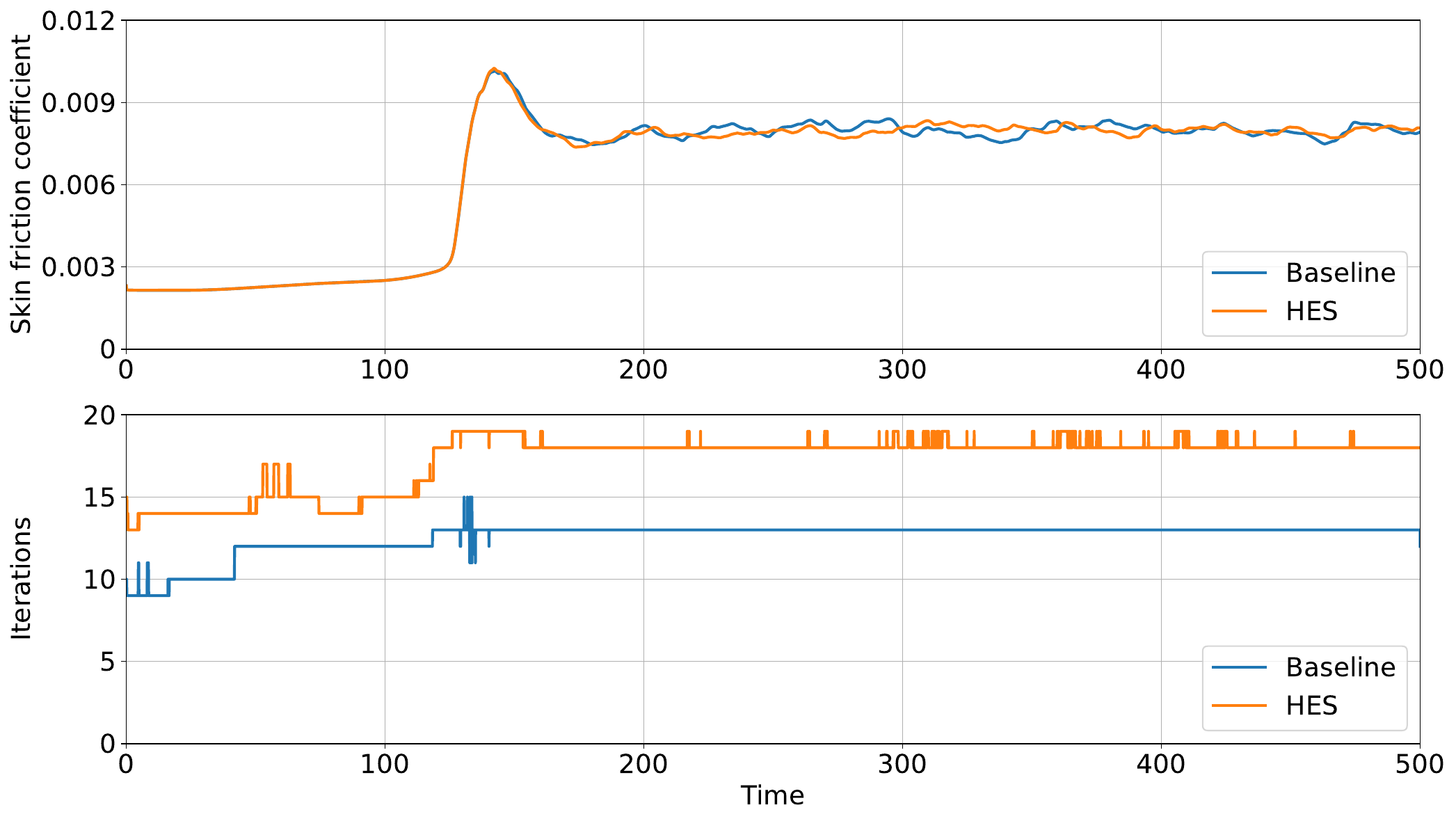}
\vspace{-0.3cm}
\caption{Evolution of the number of iterations during the simulation, plain channel flow.}
\label{fig:channel_iters}
\end{figure}

The \textit{baseline} solver configuration provides the turbulent flow calculation results in 231~min, of which 217~min are spent on solving the pressure Poisson equation and the other 14~min on calculating the spatial discretizations and performing integration in time. The distribution of the cumulative number of iterations per time step is shown in Figure~\ref{fig:channel_iters}. The number of iterations correlates with the evolution of the skin friction coefficient, and a slight increase from 10--12~to 13~iterations is observed after the transition to turbulence. 

The optimized linear solver configuration allows for reducing the cumulative SLAE solution time to 103~min and the whole simulation takes 117~min. The ratio of the \textit{baseline} and optimized SLAE solution times indicates that more than a twofold speedup can be achieved. The factual DNS speedup depends on the time integration interval (the number of time steps modeled), and it will tend to increase up to the SLAE solution speedup with increasing the time integration interval. Accounting for the extra 10~min for linear solver optimization, the overall simulation time for the test problem modeled reaches 127~min, which results in a speedup by a factor of~1.8. The evolution of the number of iterations reproduces the same behavior as for the \textit{baseline} solver configuration: the average number of iterations per time step increases from 14 to 18 when the transition from the laminar to the turbulent flow regime occurs.

%%%%%%%%%%%%%%%%%%%%%%%%%%%%%%%%%%%%%%%%

\subsection{Flow in a channel with a matrix of wall-mounted cubes}
\label{sec:channel_cube}
The second test case performs the DNS in a channel with a matrix of equally spaced cubes mounted on a channel wall. The flow is modeled in a periodic box with dimensions $4h \times 3.4h \times 4h$, where $h$ is the cube height. The modeling is performed for $\mbox{Re}_b = 3854$, defined using the bulk velocity and cube height. The overall grid size constructed is 2.32~mln. cells, and the grid stretching is applied for each wall. More details about the problem statement can be found in~\cite{KrasnopolskyCPC2018}. The simulation models $T=600$ time units, and the transition to turbulence occurs in the first 100~time units. The time integration step is set to $\tau = 4 \cdot 10^{-3}$. The tolerances for solving the linear systems are set to $10^{-8}$ and $10^{-10}$ for the pressure and pseudo-pressure, respectively. The formulated case is calculated on six compute nodes. 

\begin{figure}[t]
\centering
\includegraphics[width=11cm]{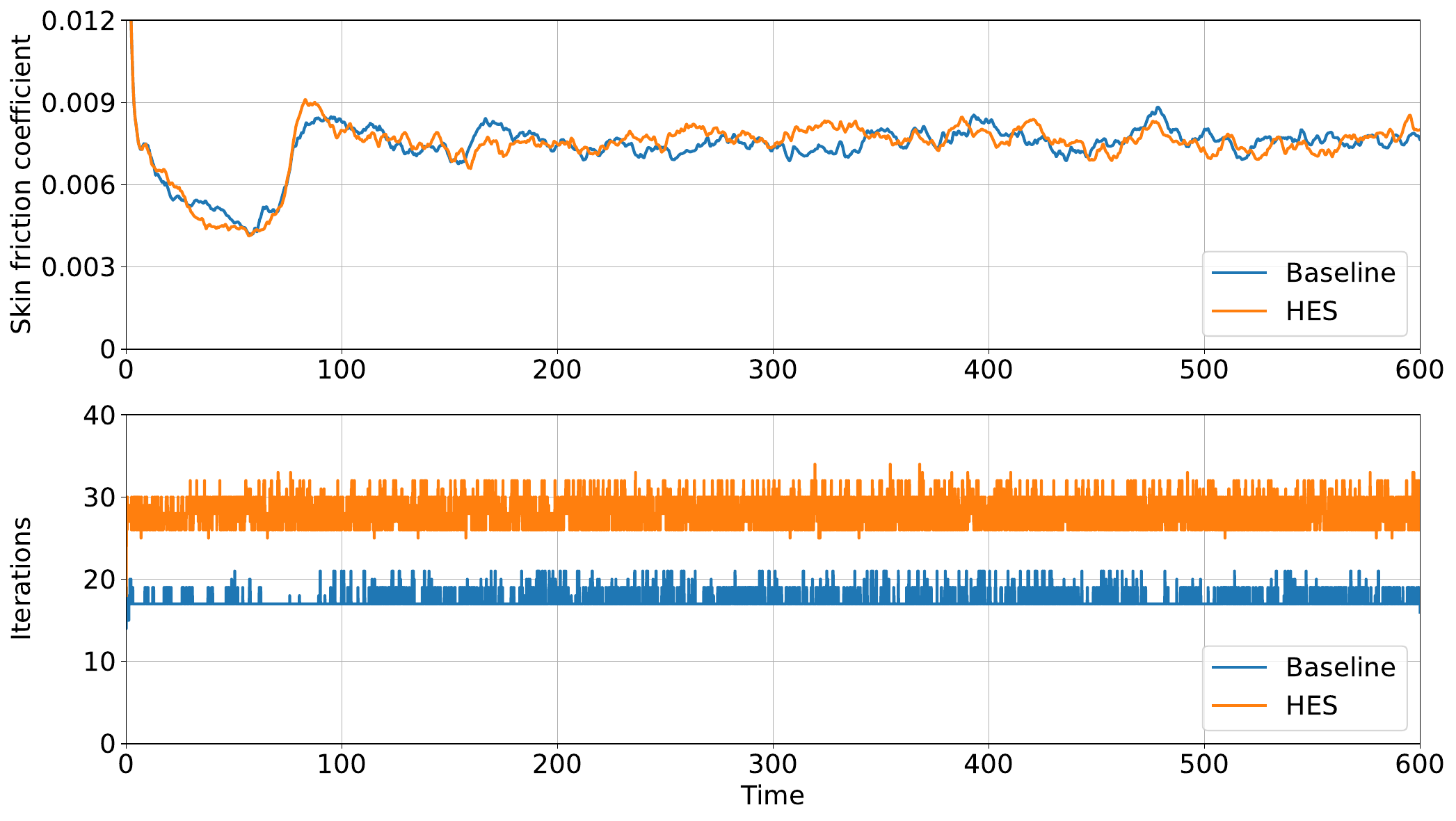}
\vspace{-0.3cm}
\caption{Evolution of the number of iterations during the simulation, flow in a channel with a matrix of wall-mounted cubes.}
\label{fig:channel_cube_iters}
\end{figure}

The simulation performed with the \textit{baseline} linear solver configuration is done in 1007~min, which includes 967~min spent on solving SLAEs and 40~min spent on the other calculations. The evolution of the cumulative number of linear solver iterations per time step shows almost constant variation during the whole run (Figure~\ref{fig:channel_cube_iters}), and the number of iterations changes in the range of 17--21. Modeling with the optimized linear solver configuration allows for a significant reduction in calculation time. The corresponding SLAE solution time is reduced to 621~min, which indicates a speedup of a factor of~1.56. The full DNS takes 662~min, and even after accounting for 12~min for SLAE solver optimization, the resulting advantage for the simulation with $T=600$ is by a factor of~1.5. The number of iterations for the optimized solver configuration varies in the range of 26--31.

%%%%%%%%%%%%%%%%%%%%%%%%%%%%%%%%%%%%%%%%

\subsection{Flow in a square duct}
\label{sec:duct}\label{Ed:summary3}
The third test problem simulates turbulent flow in a straight square duct with $\mbox{Re}_b = 2205$ ($\mbox{Re}_{\tau} = 150$). The $4\pi h \times 2h \times 2 h$ box is considered, where $h$ corresponds to the duct semi-height. The problem is solved on the non-uniform grid, consisting of $210 \times 128 \times 128$ cells. The simulation time interval of $T = 500$ time units is modeled with the time integration step $\tau = 1.75 \cdot 10^{-2}$. The SLAEs for the pressure and pseudo-pressure are solved until reaching the tolerances of $10^{-6}$ and $10^{-8}$. This test series is performed with 4~compute nodes and 3~linear solver configurations. In addition to the \textit{baseline}, the ES and HES optimized configurations are obtained and used to perform DNS. Here, due to observable variations in the ES optimization results across the several test runs, the generation size has been increased twice, and 10~soft and 10~random mutations have been performed. This modification allowed for improved convergence but led to an increase in the optimization time.

\begin{figure}[t]
\centering
\includegraphics[width=11cm]{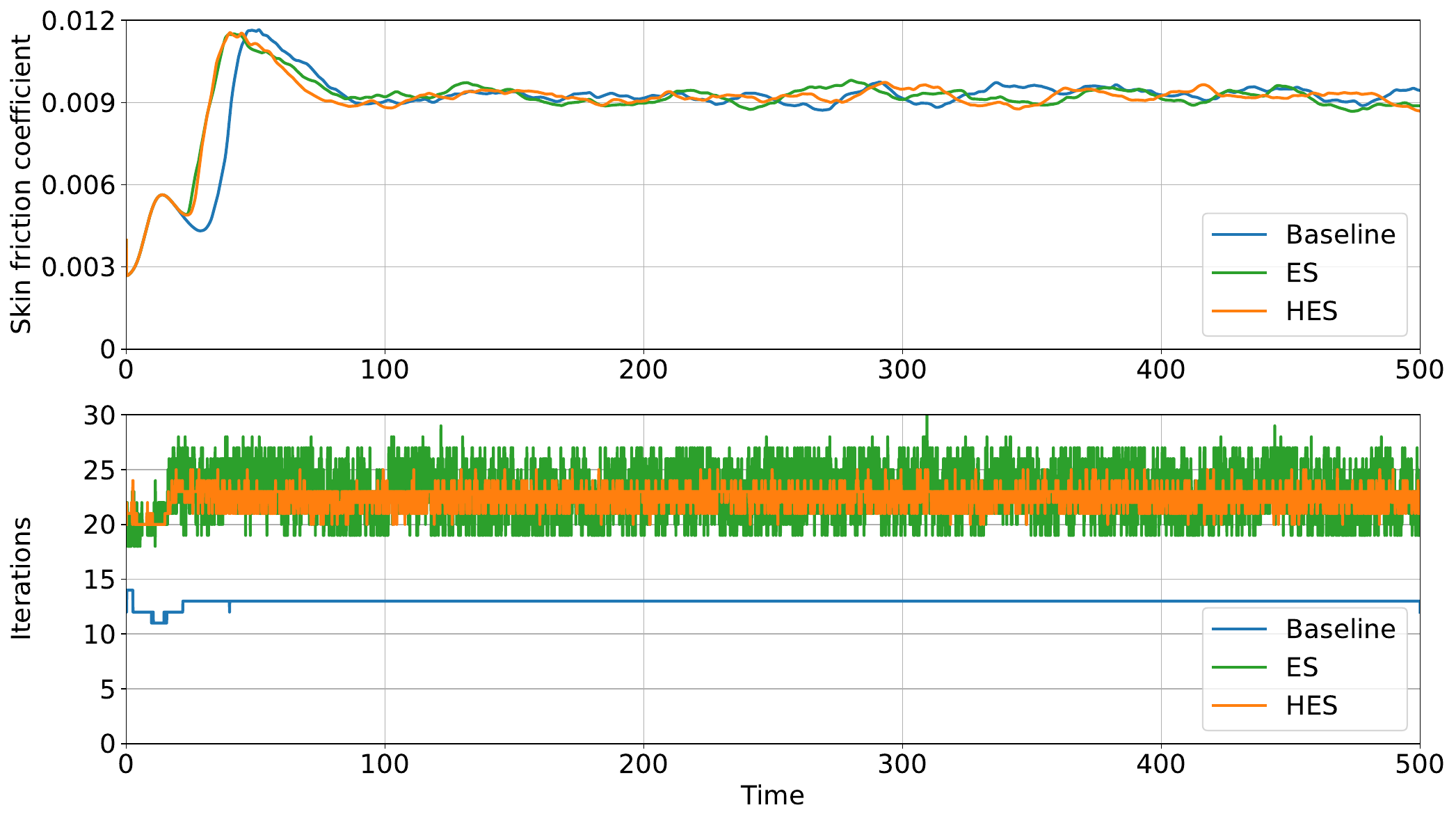}
\vspace{-0.3cm}
\caption{Evolution of the number of iterations during the simulation, flow in a square duct.}
\label{fig:duct_iters}
\end{figure}

The use of the \textit{baseline} solver configuration allows for the simulation to be completed in 283~min, spending 267~min with the linear solver and 16~min on discretizations. The typical number of iterations per time step is equal to~13 (Figure~\ref{fig:duct_iters}). The evolution strategy provides the optimized parameter configuration in 29~min; this allows for solving a sequence of SLAEs in 256~min and performing the DNS in 275~min. The cumulative number of iterations per time step for this solver configuration varies between 21 and 24. The HES optimized configuration requires 57~min for optimization (due to the twice larger number of generations compared to the ES run performed), but provides a much more efficient solver configuration capable of solving SLAEs in 189~min and calculating the flow in 207~min. This results in a decrease in the SLAE solution time by a factor of~1.41. The number of linear solver iterations until convergence varies in the range of 19--26.

The obtained results indicate that both the ES and HES produce results that outperform the \textit{baseline}. The speedup of the DNS calculation depends on the optimization time and the length of the time integration interval specified by the researcher. For the test problem considered ($T=500$), the optimization overhead for ES increases the resulting simulation time to 304~min, which nullifies the positive effect of the linear solver speedup. For the HES, the total simulation time reaches 264~min, still providing the calculation speedup by a factor of~1.07. The effect, however, will become more valuable when performing longer simulation times and/or improving the XAMG library setup phase to reduce the evolution strategy optimization time. The latter one can be significantly reduced by reworking some internal data conversion procedures in XAMG that affect the multigrid method setup phase.

%%%%%%%%%%%%%%%%%%%%%%%%%%%%%%%%%%%%%%%%

\subsection{Analysis of the DNS results}
\label{Ed:summary4}
The optimized configurations of the linear solvers obtained in the sections above require more iterations compared to the \textit{baseline} configuration. Despite the increased number of iterations, a decrease in the overall calculation time is achieved. That is a result of significantly different multigrid matrices hierarchy, constructed for the optimized linear solver configuration. The optimized configuration has a twice lower number of levels in the hierarchy, and the finest levels contain a much smaller amount of nonzero elements per matrix row. An example of the matrix hierarchies constructed for the \textit{baseline} and optimized linear solver configurations for the turbulent flow in a channel with a cube is shown in Table~\ref{tab:hierarchy_channel_cube}. These observations argue that the main factor resulting in the acceleration of the SLAE solution time is reducing the overall number of floating point operations, resulting in the algorithmic speedup but not the related software optimizations.

The calculation results presented above confirm the possibility of using the proposed hybrid evolution strategy to tune the parameters of linear solvers when modeling incompressible turbulent flow simulations. On par with the observable calculation speedup, the HES automates and, thus, dramatically simplifies the question of choosing the linear solver configuration for the specific problem of interest.

\begin{table*}[t]
\centering
\caption{Multigrid matrix hierarchies for the two linear solver configurations, flow in a channel with a matrix of wall-mounted cubes.}
\begin{tabular}{| c || c | c | c || c | c | c |}
\hline
\multirow{2}*{Level} & \multicolumn{3}{c||}{\textit{baseline}} & \multicolumn{3}{c|}{optimized} \\
\cline{2-7}
& \texttt{nrows} & \texttt{nonzeros} & \texttt{avg/row} & \texttt{nrows} & \texttt{nonzeros} & \texttt{avg/row} \\
\hline
0  & 2322432 &  15426816 & 6.6   & 2322432& 15426816 & 6.6  \\
\hline
1  & 1131273 &  17019883 & 15    & 452323 & 5625853  & 12.4 \\
\hline
2  & 344091  & 7934791   & 23.1  & 58192  & 904930   & 15.6 \\
\hline
3  & 176824  & 8068816   & 45.6  & 4443   & 68823    & 15.5 \\
\hline
4  & 92676   & 7162768   & 77.3  & 1112   & 24232    & 21.8 \\
\hline
5  & 47687   & 5658159   & 118.7 & 236    & 6618     & 28   \\
\hline
6  & 24054   & 3919982   & 163   &        &          &      \\
\hline
7  & 11682   & 2371260   & 203   &        &          &      \\
\hline
8  & 5340    & 1243284   & 232.8 &        &          &      \\
\hline
9  & 2190    & 551346    & 251.8 &        &          &      \\
\hline
10 & 827     & 206121    & 249.2 &        &          &      \\
\hline
11 & 269     & 50235     & 186.7 &        &          &      \\
\hline
\end{tabular}
\label{tab:hierarchy_channel_cube}
\end{table*}

%%%%%%%%%%%%%%%%%%%%%%%%%%%%%%%%%%%%%%%%

\section{Conclusions}

The use of advanced iterative methods for solving systems of linear algebraic equations suffers lots of difficulties due to a large number of tuning parameters affecting the corresponding methods productivity. The present paper proposes an automated parameter tuning algorithm capable of optimizing the iterative method parameters for the specific linear system of interest. The algorithm is based on the evolution strategy with the hybrid random mutation operator, involving the pre-trained neural network. The algorithm combines two distinct properties: adaptivity, realized by the evolution strategy, and a~priori knowledge, introduced by the neural network. The neural network usage scenario does not require accurate predictions for the solution times with a specific linear system and hardware platform but only some relative values. This aspect greatly simplifies the neural network training process and allows for reusing the pre-trained ones with various linear systems.

The paper provides a detailed formulation of the proposed optimization algorithm and the neural network architecture. Several methodological aspects of training the neural network accounting for the specific problem of interest are discussed. A series of numerical experiments is performed to indicate the influence of hybrid evolution strategy control parameters on the optimization results and the required calculation time. These tests demonstrate the potential of SLAE solution acceleration with the proposed evolution strategy at the cost of about 100~trial solutions. The use of the hybrid mutation operator with the neural network filtering allows for further improvement of the optimization results, a decrease in variation across the repeated launches, and a reduction in the optimization algorithm execution time.

The incompressible turbulent flow simulations are performed to indicate the DNS calculation time reduction when using the hybrid evolution strategy to tune the linear solver parameters. The speedup by a factor of 1.4-2 is achieved for solving the sequence of SLAEs for the pressure Poisson equation. This results in a speedup of up to 1.8~times for the DNS tests performed, even after accounting for the extra time spent on parameter optimization. The pre-trained networks can be successfully reused to optimize various systems on different compute platforms. Another important feature of this algorithm is the simplification of the modeling process and the automation of the manual routine work done by the researchers.

%%%%%%%%%%%%%%%%%%%%%%%%%%%%%%%%%%%%%%%%

\section*{Declaration of competing interest} 
The authors declare that they have no known competing financial interests or personal relationships
that could have appeared to influence the work reported in this paper.

%%%% 

\section*{Acknowledgments}
\label{lab:ack}
The presented work is supported by RSF grant No. 18-71-10075. The research is carried out using the equipment of the shared research facilities of HPC computing resources at Lomonosov Moscow State University and the computing resources of the federal collective usage center Complex for Simulation and Data Processing for Mega-science Facilities at the NRC ``Kurchatov Institute'', http://ckp.nrcki.ru/.

%%%%%%%%%%%%%%%%%%%%%%%%%%%%%%%%%%%%%%%%

\appendix

\section{Neural network architecture}
\label{app:NN}

The neural network is constructed in this paper to generalize and reuse the calculation results. The details of the NN architecture are summarized below.

\begin{itemize}
\item \textbf{Layers and activation functions.} The NN input layer takes the values of each parameter, and the output layer containing 1~neuron produces the solution time estimate for a given parameter vector. The three hidden dense layers, with 512, 256, and 128~neurons, respectively, follow the input layer. Hidden layers have sigmoid activation functions; the output neuron has a linear activation function. The choice of the activation functions and the size of the hidden layers have been optimized using the KerasTuner random search scheme~\cite{omalley2019kerastuner}. The dropout layers are used after each hidden layer to prevent overfitting. The dropout rate is set equal to~0.25 in accordance with the results of the grid search procedure.

\item \textbf{Loss Function.} The mean squared error~(MSE), defined as
\begin{equation}
MSE = \frac{1}{N_t} \sum_{i=1}^{N_t} \left(T^*_i - \hat{T_i} \right)^2,
\end{equation}
where $N_t$ is the size of the training dataset, $\hat{T_i}$ is the predicted value, and ${T^*_i}$ is the input data, is a commonly used loss function for regression problems~\cite{LossSurvey, LossFuncInvestigation}. The MSE loss function minimizes the difference between the predicted and actual values in the whole data range.

\item \textbf{Model optimizer.} The ADAM (Adaptive Moment Estimation) optimizer is used for training the NN model. The initial learning rate is set to $10^{-3}$. The number of epochs as a function of the training dataset size is selected in accordance with the following empirical formula:
\begin{equation}
N_{epochs} = \frac{N_t}{50}+50.
\end{equation}
This estimate ensures stable training of the constructed model in the range of used dataset sizes ($N_t \sim 10^4-10^5$) and the number of parameters ($N \sim 20$) for the test problems performed.

\item \textbf{Cross-validation technique.} The data cross-validation technique is used to provide an effective use of the input dataset and a representative assessment of the accuracy of the model. For the proposed model, a single hold-out random subsampling method is applied~\cite{crossval}. The input dataset is split into the training and validation parts. The corresponding method applies a random permutation of the samples in the input dataset, thus mixing the data between the training and validation sets at the subsequent training sessions.
\end{itemize}

%%%%%%%%%%%%%%%%%%%%%%%%%%%%

\section{List of optimized parameters}
\label{app:list}
The present paper deals with the optimization of four groups of parameters containing 7, 13, 15, and 21~numerical method parameters. These parameters correspond to the BiCGStab iterative method and the classical algebraic multigrid preconditioner. The Chebyshev polynomial method is chosen as both the pre- and post-smoother (with their own method parameters), and the direct solver is used at the coarsest level; these basic methods are fixed and do not change during the optimization. The full list of the XAMG library parameters targeted for optimization and their equivalents in \textit{hypre} are summarized in Table~\ref{tab:params}. The specific parameters included in the four optimization sets are shown in Table~\ref{tab:4sets}.

Table~\ref{tab:params} also contains the \textit{baseline} and \textit{default} solver configurations used in this paper. The \textit{baseline} linear solver configuration used for the tests is a universal and robust one formulated as a result of some modeling experience. The \textit{default} configuration is constructed based on the recommendations provided in the \textit{hypre} library documentation. A slightly better variant can probably be proposed for the problems considered in the paper as a result of thorough tuning. This may slightly change the reference values and speedup numbers, but it does not affect the key results discussed in the paper.

\begin{table*}[!hbt]
\caption{The list of parameters assigned for optimization and the predefined parameter configurations.}
\begin{center}
\begin{tabular}{|l | l | l | c | c |}
\hline
XAMG & \textit{hypre} & Allowed values & \textit{baseline} & \textit{default} \\
\hline
\multicolumn{5}{|l|} {\hspace{0.5cm}Preconditioner: \texttt{Multigrid}} \\
\hline
\texttt{max\_iters}  & MaxIter & 1, 2, 3 & 1 & 1 \\
\hline
\texttt{mg\_cycle}  & CycleType & $V, W, F^*$ & $V$ & $V$ \\
\hline
\texttt{mg\_max\_row\_sum}  & MaxRowSum & $\varepsilon^{**}$, 0.05, \ldots, 1 & 1 & 0.9 \\
\hline
\texttt{mg\_nonGalerkin\_tol}  & NonGalerkinTol & 0, 0.05, \ldots, 1 & 0 & 0 \\
\hline
\texttt{mg\_coarse\_matrix\_size}  & MaxCoarseSize & 50, 100, \ldots, 500 & 500 & 9 \\
\hline
\texttt{mg\_num\_paths}  & NumPaths & 1-4 & 3 & 1 \\
\hline
\texttt{mg\_coarsening\_type}  & CoarsenType & 0, 3, 6, 8 & 6 & 6 \\
\hline
\texttt{mg\_interpolation\_type}  & InterpType & 0, 3-9, 14 & 0 & 0 \\
\hline
\texttt{mg\_strength\_threshold}  & StrongThreshold & 0, 0.1, \ldots, 0.9 & 0.25 & 0.5 \\
\hline
\texttt{mg\_trunc\_factor}  & TruncFactor & 0, 0.1, \ldots, 0.9 & 0.3 & 0 \\
\hline
\texttt{mg\_Pmax\_elements}  & PMaxElmts & 0-10 & 4 & 4 \\
\hline
\texttt{mg\_agg\_num\_levels}  & AggNumLevels & 0-10 & 2 & 0 \\
\hline
\texttt{mg\_agg\_interpolation\_type}  & AggInterpType & 1-4 & 4 & 1 \\
\hline
\texttt{mg\_agg\_trunc\_factor}  & AggTruncFactor & 0, 0.1, \ldots, 0.9 & 0.3 & 0 \\
\hline
\texttt{mg\_agg\_P12\_trunc\_factor}  & AggP12TruncFactor & 0, 0.05, \ldots, 1 & 0 & 0 \\
\hline
\texttt{mg\_agg\_Pmax\_elements}  & AggPMaxElmts & 1-10 & 4 & 4 \\
\hline
\texttt{mg\_agg\_P12max\_elements}  & AggP12MaxElmts & 1-10 & 4 & 4 \\
\hline
\multicolumn{5}{|l|} {\hspace{0.5cm}Smoother: \texttt{Chebyshev}} \\
\hline
\texttt{polynomial\_order}  & ChebyOrder & 1-4 & 2 & 2 \\
\hline
\texttt{spectrum\_fraction}  & ChebyFraction & $\varepsilon$, 0.1, \ldots, 0.9 & 0.3 & 0.3 \\
\hline
\end{tabular}
\end{center}
 * not supported by \textit{hypre}
\\ ** $\varepsilon$ is set equal to $10^{-10}$
\label{tab:params}
\end{table*}

\begin{table*}[!hbt]
\caption{Lists of parameters included in the optimization sets.}
\begin{center}
\begin{tabular}{| l | c | c | c | c |}
\hline
\multirow{2}*{Parameter} & \multicolumn{4}{|c|}{Optimized set} \\
\cline{2-5}
 & 7 parameters & 13 parameters & 15 parameters & 21 parameters \\
\hline
% \multicolumn{4}{|l|} {\hspace{0.5cm}Preconditioner: \texttt{Multigrid}} \\
% \hline
\texttt{max\_iters}  & 1 & 1 & * & * \\
\hline
\texttt{mg\_cycle}  & $V$ & $V$ & * & * \\
\hline
\texttt{mg\_max\_row\_sum}  & 0.9 & * & 0.9 & * \\
\hline
\texttt{mg\_nonGalerkin\_tol}  & 0.0 & * & 0.0 & * \\
\hline
\texttt{mg\_coarse\_matrix\_size}  & 100 & * & 100 & * \\
\hline
\texttt{mg\_num\_paths}  & * & * & * & * \\
\hline
\texttt{mg\_coarsening\_type}  & * & * & * & * \\
\hline
\texttt{mg\_interpolation\_type}  & * & * & * & * \\
\hline
\texttt{mg\_strength\_threshold}  & 0.5 & * & 0.5 & * \\
\hline
\texttt{mg\_trunc\_factor}  & 0.25 & * & 0.25 & * \\
\hline
\texttt{mg\_Pmax\_elements}  & 4 & 4 & * & * \\
\hline
\texttt{mg\_agg\_num\_levels}  & * & * & * & * \\
\hline
\texttt{mg\_agg\_interpolation\_type}  & * & * & * & * \\
\hline
\texttt{mg\_agg\_trunc\_factor}  & 0.25 & * & 0.25 & * \\
\hline
\texttt{mg\_agg\_P12\_trunc\_factor}  & 0 & 0 & * & * \\
\hline
\texttt{mg\_agg\_Pmax\_elements}  & 4 & 4 & * & * \\
\hline
\texttt{mg\_agg\_P12max\_elements}  & 4 & 4 & * & * \\
\hline
% \multicolumn{4}{|l|} {\hspace{0.5cm}Smoother: \texttt{Chebyshev}} \\
% \hline
\texttt{polynomial\_order}$^{**}$  & * & * & * & * \\
\hline
\texttt{spectrum\_fraction}$^{**}$  & 0.3 & 0.3 & * & * \\
\hhline{|=|=|=|=|=|}
\textit{Total combinations:} & $10^6$ & $4\cdot10^{12}$ & $2 \cdot 10^{12}$ & $7\cdot10^{18}$ \\
\hline
\end{tabular}
\end{center}
 * parameter is assigned for optimization \\ 
 ** independent parameters for pre- and post-smoother
\label{tab:4sets}
\end{table*}

%%%% 

\bibliographystyle{model1-num-names}
\bibliography{optim}

%%%%%%%%%%%%%%%%%%%%%%%%%%%%%%%%%%%%%%%%%%%%%%%%%%%%%%
%%%%%%%%%%%%%%%%%%%%%%%%%%%%%%%%%%%%%%%%%%%%%%%%%%%%%%
%%%%%%%%%%%%%%%%%%%%%%%%%%%%%%%%%%%%%%%%%%%%%%%%%%%%%%
%%%%%%%%%%%%%%%%%%%%%%%%%%%%%%%%%%%%%%%%%%%%%%%%%%%%%%
%%%%%%%%%%%%%%%%%%%%%%%%%%%%%%%%%%%%%%%%%%%%%%%%%%%%%%

\iffalse
\newpage
\begin{nolinenumbers}

\setcounter{page}{1}

\include{reply_to_reviewers}

\end{nolinenumbers}
\fi

\end{document}